\input ./style/arxiv-general.cfg
\documentclass[aap,MSNbibl,dvips]{arximspdf}
\makeatletter
   \@ifpackageloaded{graphicx}{}{\usepackage{graphicx}}
\makeatother

\doi{10.1214/15-AAP1137}
\volume{26}
\issue{4}
\pubyear{2016}
\firstpage{1947}
\lastpage{1995}
\docsubty{FLA}

\makeatletter
\def\cal{\mathcal}

\newcommand{\eqref}[1]{(\ref{#1})}

\newcommand{\ME}{{\mathbb E}}

\newtheorem{proposition}{Proposition}[section]
\newtheorem{corollary}[proposition]{Corollary}
\newtheorem{lemma}[proposition]{Lemma}

\newproclaim{remark}[proposition]{Remark}
\newproclaim{assumption}[proposition]{Assumption}

\newproclaim{condition}[proposition]{Condition}
\newproclaim{example}[proposition]{Example}

\newproclaim{definition}[proposition]{Definition}

\makeatother

\begin{document}
\begin{frontmatter}

\title{Asymptotically optimal priority policies
for indexable and nonindexable restless bandits}
\runtitle{Asymptotic optimal control of restless bandits}

\begin{aug}
\author[A]{\fnms{I.~M.}~\snm{Verloop}\corref{}\ead
[label=e1]{verloop@irit.fr}}
\runauthor{I.~M. Verloop}
\affiliation{CNRS, IRIT, Universit\'e de Toulouse and INP}
\address[A]{IRIT\\
2 rue C. Camichel\\
F-31071 Toulouse\\
France\\
\printead{e1}}
\end{aug}

%
\received{\smonth{6} \syear{2014}}
%
\revised{\smonth{8} \syear{2015}}

%
\begin{abstract}
We study the asymptotic optimal control of multi-class restless
bandits. A restless bandit is a controllable stochastic process whose
state evolution depends on whether or not the bandit is made active.
Since finding the optimal control is typically intractable, we propose
a class of priority policies that are proved to be asymptotically
optimal under a global attractor property and a technical condition.
We consider both a fixed population of bandits as well as a dynamic
population where bandits can depart and arrive.
As an example of a dynamic population of bandits, we analyze a
multi-class $\mathit{M/M/S+M}$ queue for which we show asymptotic
optimality of an index policy.

We combine fluid-scaling techniques with linear programming results to
prove that when bandits are indexable, Whittle's index policy is
included in our class of priority policies. We thereby generalize a
result of Weber and Weiss [\textit{J. Appl. Probab.} \textbf{27} (1990) 637--648]
 about asymptotic optimality of
Whittle's index policy to settings with (i) several classes of bandits,
(ii) arrivals of new bandits and (iii) multiple actions.

Indexability of the bandits is not required for our results to hold.
For nonindexable bandits, we describe how to select priority policies
from the class of asymptotically optimal policies and present numerical
evidence that, outside the asymptotic regime, the performance of our
proposed priority policies is nearly optimal.
\end{abstract}

%

%
\end{frontmatter}

\section{Introduction}\label{sec1}

Multi-armed bandit problems are concerned with the optimal dynamic
activation of several competing \emph{bandits}, taking into account
that at each moment in time $\alpha$ bandits can be made active.
A bandit is a controllable stochastic process whose state evolution
depends on whether or not the bandit is made active.
The aim is to find a control that determines at each decision epoch
which bandits to activate in order to minimize the overall cost
associated to the states the bandits are in.
In the by now classical multi-armed bandit model, \cite{Gittins89}, it
is assumed that only active bandits can change state. In \cite
{Whittle88}, Whittle introduced the so-called restless bandits, where a
bandit can also change its state while being passive (i.e., not active),
possibly according to a different law from the one that applies when
it is active.
The multi-armed restless bandit problem is a stochastic optimization
problem that has gained popularity due to its multiple applications in,
for example, sequential selection trials in medicine, sensor
management, manufacturing systems, queueing and communication networks,
control theory, economics, etc.
We refer to \cite{GGW11,MT07,Whittle96} for further references,
applications, and possible extensions that have been studied in the literature.

In 1979, Gittins \cite{Gittins79} introduced index-based policies for
the nonrestless bandit problem.
He associated to each bandit an index, which is a function of the state
of the bandit, and defined the policy that
activates $\alpha$ bandits with currently the largest indices.
This policy is known as the Gittins index policy.
It was first proved by Gittins that this policy is optimal in the case
$\alpha=1$ \cite{Gittins79} for the time-average and discounted cost criteria.
In the presence of restless bandits, finding an optimal control is
typically intractable. In 1988, Whittle \cite{Whittle88} proposed
therefore to solve a relaxed optimization problem where the constraint
of having at most $\alpha$ bandits active at a time is relaxed to a
time-average or discounted version of the constraint.
In addition, he defined the so-called indexability property, which
requires to establish that as one increases the Lagrange multiplier of
the relaxed optimization problem, the collection of states in which the
optimal action is passive increases. Under this property, Whittle
showed that an optimal solution to the relaxed optimization problem can
be described by index values. The latter, in turn, provide a heuristic
for the original restless bandit problem, which is referred to as
Whittle's index policy in the literature. It reduces to Gittins index
policy when passive bandits are static (the nonrestless case).
Whittle's index policy is in general not an optimal solution for the
original problem.
In \cite{WW90}, Weber and Weiss proved Whittle's index policy to be
asymptotically optimal.


In this paper, we study the asymptotic optimal control of a general
multi-class restless bandit problem. We consider both a fixed
population of bandits as well as a dynamic scenario where bandits can
arrive and depart from the system.
The asymptotic regime is obtained by letting the number of bandits that
can be simultaneously made active grow proportionally with the
population of bandits.

In one of our main results, we derive a set of priority policies that
are asymptotically optimal when certain technical conditions are satisfied.
In another main result, we then prove that if the bandits are
indexable, Whittle's index policy is contained in our set of priority
policies. We thereby generalize the asymptotic optimality result of
Weber and Weiss \cite{WW90} to settings with (i) several classes of
bandits, and (ii) arrivals of new bandits.
Another extension presented in the paper is the possibility of choosing
among multiple actions per bandit. This is referred to as
``super-process'' in the literature \cite{GGW11}.
Throughout the paper, we discuss how our asymptotic optimality results
extend to that scenario.

Efficient control of \emph{nonindexable} restless bandits has so far
received little attention in the literature. Nonindexable settings can
however arise in problems of practical interest; see, for
example, \cite
{HG11} in the context of a make-to-stock system.
The definition of our set of priority policies does not rely on
indexability of the system, and hence, provides asymptotically optimal
heuristics for nonindexable settings.
We describe how to select priority policies from this set and present
numerical evidence that, outside the asymptotic regime, the performance
of our proposed priority policies is nearly optimal.


The asymptotic optimality results obtained in this paper hold under
certain technical conditions. For a fixed population of bandits, these
conditions reduce to a differential equation having a global attractor,
which coincides with the condition as needed by Weber and Weiss \cite{WW90}.
For a dynamic population of bandits, additional technical conditions
are needed due to the infinite state space.
To illustrate the applicability of the results, we present a large
class of restless bandit problems for which we show the additional
technical conditions to hold. This class is characterized by the fact
that a bandit that is kept passive will eventually leave the system.
This can represent many practical situations such as impatient
customers, companies that go bankrupt, perishable items, etc.
We then present a multi-class $\mathit{M/M/S+M}$ queue, which is a very
particular example of the above described class. We describe a priority
policy that satisfies the global attractor property, and hence
asymptotic optimality follows.

In this paper, we consider a generalization of the standard restless
bandit formulation: Instead of having at each moment in time exactly
$\alpha$ bandits active, we allow strictly less than $\alpha$ bandits
to be active at a time. We handle this by introducing so-called dummy
bandits. In particular, we show that it is asymptotically optimal to
activate those bandits having currently the largest, \emph{but strictly
positive}, Whittle's indices. Hence, whenever a bandit is in a state
having a negative Whittle's index, this bandit will never be activated.

Our proof technique relies on a combination of fluid-scaling techniques
and linear programming results: First, we describe the fluid dynamics
of the restless bandit problem, taking only into account the average
behavior of the original stochastic system.
The optimal equilibrium points of the fluid dynamics are described by
an LP problem. We prove that the optimal value of the LP provides a
lower bound on the cost of the original stochastic system.
The optimal fluid equilibrium point is then used to describe priority
policies for the original system whose fluid-scaled cost coincides with
the lower bound, and are hence referred to as asymptotically optimal policies.
In order to prove that Whittle's index policy is one of these
asymptotically optimal policies, we then reformulate the relaxed
optimization problem into an LP problem. An optimal solution of this LP
problem is proved to coincide with that of the LP problem corresponding
to the fluid problem as described above.
This is a different proof approach than that taken in \cite{WW90} and
allows to include arrivals of bandits to the system, whereas the
approach of \cite{WW90} does not.




To summarize, the main contributions of this paper are the following:
\begin{itemize}
\item For a multi-class restless bandit problem (possibly
nonindexable) with either a fixed or dynamic population of bandits, we
determine a set of priority policies that are asymptotically optimal if
the corresponding ODE has a global attractor and certain technical
conditions (Condition \ref{cond:twee}) are satisfied (Proposition~\ref
{45454}).
\item We show that Condition \ref{cond:twee} is satisfied for a large
class of restless bandit problems. In particular, for a fixed
population of bandits under a unichain assumption and for a dynamic
population when passive bandits will eventually leave the system
(Proposition~\ref{35353}).
%
\item In the case the bandits are indexable, we show that Whittle's
index policy is inside our set of priority policies, both for a fixed
population of bandits (Proposition~\ref{56565}) and for a
dynamic population of bandits (Proposition~\ref{75757}).
%
\item For nonindexable bandits, we describe how to select priority
policies from the class of asymptotically optimal policies
(Section~\ref{sec:nonindex_serie}) and for a particular example we numerically show
that outside the asymptotic regime their sub-optimality gap is very
small (Section~\ref{sec:nonindex_example}).
\end{itemize}

The remainder of the paper is organized as follows. In Section~\ref{sec:relatedwork}, we give an overview of related work and in
Section~\ref{sec:model} we define the multi-class restless bandit problem.
In Section~\ref{sec:fluid}, we define our set of priority policies and
state the asymptotic optimality result, both for a fixed population as
well as for a dynamic population of bandits.
In Section~\ref{sec:whittle}, we define Whittle's index policy and
prove it to be asymptotically optimal.
In Section~\ref{sec:examples}, we discuss the global attractor property
required in order to prove the asymptotic optimality result.
In Section~\ref{sec:abopt}, we present the $\mathit{M/M/S+M}$ queue as an
example of an indexable restless bandit and derive a robust priority
policy that is asymptotically optimal.
Section~\ref{sec:nonindex} focuses on the selection of asymptotically
optimal priority policies for nonindexable bandits and numerically
evaluates the performance.

\section{Related work}
\label{sec:relatedwork}

For the nonrestless bandit problem, optimality of Gittins index policy
has been proved in \cite{Gittins79}, for the case $\alpha=1$ and a
time-average or discounted cost criteria.
In \cite{Weiss88,Whittle81}, the optimality result was extended to a
dynamic population of bandits where new bandits may arrive over time,
for example, Poisson arrivals or Bernouilli arrivals.
For $\alpha>1$, the optimality results do not necessarily go through.
In \cite{PT99}, sufficient conditions on the reward processes were
given in order to guarantee optimality of the Gittins policy for the
discounted cost criterion, when $\alpha>1$.

For the restless bandit problem, the authors of \cite{GHK11} have
extended Whittle's index heuristic to the setting where each restless
bandit may choose from multiple actions, that is, representing a
divisible resource to a collection of bandits.
Over the years, Whittle's index policy has been extensively applied and
numerically evaluated in various application areas such as wireless
downlink scheduling \cite{AEJ10,OES12}, systems with delayed state
observation \cite{EM04}, broadcast systems \cite{RBCK08}, multi-channel
access models \cite{ALJZK09,LZ10}, stochastic scheduling
problems \cite
{AGNK03,GM02,NM07} and scheduling in the presence of impatient
customers \cite{AJN10,GKO09,LAV14,N07}.

As opposed to Gittins policy, Whittle's index policy is in general not
an optimal solution for the original problem. For a fixed population of
bandits, optimality has been proved though for certain settings. For
example, in \cite{ALJZK09,LZ10} this has been proved for a restless
bandit problem modeling a multi-channel access system.
For a general restless bandit model, in \cite{Jacko11} Whittle's index
policy has been shown to be optimal for $\alpha=1$ when (i) there is
one dominant bandit or when (ii) all bandits immediately reinitialize
when made passive.
Other results on optimality of Whittle's index policy for a fixed
population of bandits exist for \emph{asymptotic regimes}. In \cite
{Whittle88}, Whittle conjectured that Whittle's index policy is nearly
optimal as the number of bandits that can be simultaneously made active
grows proportionally with the total number of bandits in the system.
In the case of \emph{symmetric} bandits, that is, all bandits are
governed by the same transition rules, this conjecture was proved by
Weber and Weiss \cite{WW90} assuming that the differential equation
describing the fluid approximation of the system has a global attractor.
They further presented an example for which the conjecture does not hold.
In \cite{HG15}, the approaches of \cite{WW90} were set forth and
extended to problems for which multiple activation levels are permitted
at any bandit.
Another recent result on asymptotic optimality can be found in \cite
{OES12} where the authors considered a specific model, as studied
in \cite{LZ10}, with two classes of bandits. They proved asymptotic
optimality of Whittle's index policy under a recurrence condition. The
latter condition replaces the global attractor condition needed
in \cite
{WW90} and was numerically verified to hold for their model.

For a dynamic population of restless bandits, that is, when new bandits
can arrive to the system, there exist few papers on the performance of
index policies. We refer to \cite{AEJ10,AEJV13} and \cite{AJN10} where
this has been studied in the context of wireless downlink channels and
queues with impatient customers, respectively.
In particular, in \cite{AEJ10,AJN10}, Whittle's index policy was
obtained under the discounted cost criterion and numerically shown to
perform well.
In \cite{AEJV13}, it was shown that this heuristic is in fact maximum
stable and asymptotically fluid optimal.
We note that the asymptotic regime studied in \cite{AEJV13} is
different than the one as proposed by Whittle \cite{WW90}.
More precisely, in \cite{AEJV13} at most one bandit can be made active
at a time (the fluid scaling is obtained by scaling both space and
time), while in \cite{WW90} (as well as in this paper) the number of
active bandits scales.

Arrivals of new ``entities'' to the system can also be modelled by a
fixed population of restless bandits. In that case, a bandit represents
a certain type of entities, and the state of a bandit represents the
number of this type of entities that are present in the system. Hence,
a new arrival of an entity will change the state of the bandit. In the
context of queueing systems this has been studied, for example,
in \cite
{AGNK03,GKO09,LAV14}.
A Whittle's index obtained from the relaxation of this problem
formulation can depend both on the arrival characteristics and on the
state, that is, the number of entities present in the system.
This in contrast to the dynamic population formulation of the problem,
as discussed in the previous paragraph, where the index will be
independent of the arrival characteristics or number of bandits present.
Asymptotic optimality results for a fixed population of bandits
modeling arrivals of new ``entities'' have been obtained in, for
example, \cite{GKO09} where Whittle's index was shown to be optimal
both in the light-traffic and the heavy-traffic limit.

This paper presents heuristics for nonindexable bandits that are
asymptotically optimal. Other heuristics proposed for nonindexable
problems can be found in \cite{BNM00,HG11}. In \cite{BNM00}, the
primal--dual index heuristic was defined and proved to have a
sub-optimality guarantee.
In Remark~\ref{rem:NM}, we will see that an adapted version of the
primal--dual index heuristic is included in the set of priority policies
for which we obtain asymptotic optimality results.
Using fair charges, the authors of \cite{HG11} proposed heuristics for
nonindexable bandits in the context of a make-to-stock system.
Numerically, the heuristic was shown to perform well. It can be checked
that their heuristic is not inside the set of priority policies for
which we show asymptotic optimality results.

We conclude this related work section with a discussion on the use of
LP techniques in the context of restless bandits.
An LP-based proof approach was previously used in, for example, \cite
{BNM00,NM01,NM07e}.
In \cite{NM01,NM07e}, it allowed to characterize and compute
indexability of restless bandits.
In \cite{BNM00}, a set of LP relaxations was presented, providing
performance bounds for the restless bandit problem under the
discounted-cost criterion.

\section{Model description}
\label{sec:model}

We consider a multi-class restless bandit problem in continuous time.
There are $K$ classes of bandits.
New class-$k$ bandits arrive according to a Poisson process with
arrival rate $\lambda_k\geq0$, $k=1,\ldots, K$.
At any moment in time, a class-$k$ bandit is in a certain state $j\in\{
1,2,\ldots, J_k\}$, with $J_k<\infty$.
When a class-$k$ bandit arrives, with probability $p_{k}(j)$ this
bandit starts in state $j\in\{1,\ldots, J_k\}$.

At any moment in time, a bandit can either be kept passive or active,
denoted by $a=0$ and $a=1$, respectively.
When action $a$ is performed on a class-$k$ bandit in state $i$,
$i=1,\ldots, J_k$, it makes a transition to state $j$ after an
exponentially distributed amount of time with rate $q_k(j|i,a)$,
$j=0,1,\ldots,J_k$, $j\neq i$. Here, $j=0$ is interpreted as a
departure of the bandit from the system.
We further define $q_k(j|j,a):=-\sum_{i=0, i\neq j}^{J_k}q_k(i|j,a)$.
The fact that the state of a bandit might evolve even under the passive
action explains the term of a \emph{restless} bandit.


Decision epochs are defined as the moments when an event takes place,
that is, an arrival of a new bandit, a change in the state of a bandit,
or a departure of a bandit.
A \emph{policy} determines at each decision epoch which bandits are
made active, with the restriction that at most $\alpha$ bandits can be
made active at a time.
This is a generalization of the standard restless bandit formulation
where at each moment in time \emph{exactly} $\alpha$ bandits need to be
activated, as will be explained in Remark~\ref{rem:standard_exact}.

Throughout this paper, we will consider both a fixed population of
bandits and a dynamic population of bandits:\vadjust{\goodbreak}
%
\begin{itemize}
\item\emph{Fixed population}: In this case, there are no new arrivals
of bandits, that is, $\lambda_k=0$, for all $k=1,\ldots,K$, and there
are no departures, that is, $q_k(0|j,a)=0$, for all $j,k,a$.


\item\emph{Dynamic population}: In this case, there are new arrivals
of bandits, that is, $\lambda_k>0$, for all $k=1,\ldots, K$, and each
bandit can depart from the system, that is, for each class $k$ there is
at least one state $j$ and one action $a$ such that $q_k(0|j,a)>0$.
\end{itemize}


For a given policy $\pi$, we define $X^{\pi}(t):=(X_{j,k}^{\pi,a}(t);
k=1,\ldots, K, j=1,\ldots, \break  J_k, a=0, 1)$, with $X_{j,k}^{\pi, a}(t)$
the number of class-$k$ bandits at time $t$ that are in state $j$ and
on which action $a$ is performed.
We further denote by $X^\pi_{j,k}(t):=\sum_{a=0}^{1} X_{j,k}^{\pi
,a}(t)$ the total number of class-$k$ bandits in state $j$ and $X_k^\pi
(t):=\sum_{j=1}^{J_k} X_{j,k}^\pi(t)$ the total number of class-$k$ bandits.



Our performance criteria are stability and long-run average holding cost.

\textit{Stability}.
For a given policy $\pi$, we will call the system \emph{stable} if the
process $X^\pi(t)$ has a unique invariant probability distribution.
We further use the following weaker notions of stability: a policy is
\emph{rate-stable} if\break $\lim_{t\to\infty} \sum_{j,k} \frac
{X_{j,k}^\pi
(t)}{t}=0$ almost surely and \emph{mean rate-stable} if\break $\lim_{t\to
\infty} \sum_{j,k}\frac{ \ME(X_{j,k}^\pi(t))}{t}=0$.
For a fixed population of bandits the state space is finite, hence the
process $X^\pi(t)$ being unichain is a sufficient condition for
stability of the policy $\pi$. In the case of a dynamic population of
bandits, the stability condition is more involved. Whether or not the
system is stable can depend strongly on the employed policy. In
Section~\ref{sec:fluid}, we will state necessary stability conditions
for the dynamic restless bandit problem.

\textit{Long-run average holding cost}.
Besides stability, another important performance measure is the
average holding cost.
We denote by $C_{k}(j,a)\in\mathbb{R}$, $j=1,\ldots,J_k$, the holding
cost per unit of time for having a class-$k$ customer in state $j$
under action $a$.
We note that $C_{k}(j,a)$ can be negative, that is, representing a reward.
We further introduce the following notation for long-run average
holding costs under policy $\pi$ and initial state $x:=(x_{j,k};
k=1,\ldots, K, j=1,\ldots, J_k)$:
\[
V_-^\pi(x):= \liminf_{T\to\infty} \frac{1}{T}
\ME_x \Biggl(\int_{t=0}^T \sum
_{k=1}^K \sum
_{j=1}^{J_k} \sum_{a=0}^{1}
C_{k}(j,a) X_{j,k}^{\pi
,a}(t)\,\mathrm{d}t \Biggr)
\]
and
\[
V_+^\pi(x) := \limsup_{T\to\infty} \frac{1}{T}
\ME_x \Biggl(\int_{t=0}^T \sum
_{k=1}^K \sum
_{j=1}^{J_k} \sum_{a=0}^{1}
C_{k}(j,a) X_{j,k}^{\pi
,a}(t)\,\mathrm{d}t \Biggr).
\]
If $V_-^\pi(x)=V_+^\pi(x)$, for all $x$, then we define $V^\pi
(x):=V_+^\pi(x)$.
We focus on Markovian policies, which base their decisions on the
current state and time. Our objective is to find a policy $\pi^*$ that
is average optimal, that is,
%
\begin{equation}
\label{eq:avopt} V_+^{\pi^*}(x) \leq V_-^{\pi}(x) \qquad\mbox{for all
$x$ and for all policies $\pi$,}
\end{equation}
%
under the constraint that at any moment in time at most $\alpha$
bandits can be made active, that is,
%
\begin{equation}
\label{eq:M} \sum_{k=1}^K \sum
_{j=1}^{J_k} X_{j,k}^{\pi,1}(t) \leq
\alpha \qquad\mbox {for all } t.
\end{equation}


\begin{remark}
\label{rem:standard_exact}
The standard formulation for the restless bandit problem with a fixed
population of bandits is to make \emph{exactly} $\alpha$ bandits active
at any moment in time. This setting can be retrieved from our
formulation by replacing $C_{k}(j,0)$ with $C_{k}(j,0)+C$, for all
$j,k$, where $C$ represents an additional cost of having a passive
bandit. The average additional cost for having passive bandits in the
system is equal to $(N-A)C$, with $N$ the total number of bandits in
the system and $A$ the average number of active bandits in the system.
When $C$ is large enough, an optimal policy will set $A$ maximal, that
is $A=\alpha$. Hence, we retrieve the standard formulation.
\end{remark}

\begin{remark}[(Multi actions)]
\label{rem:multi}
In the model description, we assumed there are only two possible
actions per bandit: $a=0$ (passive bandit) and $a=1$ (active bandit). A
natural generalization is to consider multiple actions per bandit, that
is, for a class-$k$ bandit in state $j$ the scheduler can chose from
any action $a\in\{0, \ldots, A_{k}(j)\}$ and at most $\alpha$ bandits
can be made active at a time, that is, $\sum_{k=1}^K \sum_{j=1}^{J_{k}}\sum_{a=1}^{A_{k}(j)} X_{j,k}^a(t)\leq\alpha$. This is
referred to as ``super-process'' in the literature \cite{GGW11}.
For the \emph{nonrestless} bandit problem with $\alpha=1$, an index
policy is known to be optimal in the case each state has a dominant
action, that is, if an optimal policy selects a class-$k$ bandit in
state $j$ to be made active, it always chooses the same action
$a_{k}(j)$, with $a_{k}(j)\in\{1, \ldots, A_k(j)\}$. A less strict
condition is given in \cite{GGW11}, Condition D.

In this paper, we focus on the setting $A_{k}(j)=1$, however, all
results obtained will go through in the multi-action context when the
definition of the policies are modified accordingly; see Remarks \ref
{rem:super1} and~\ref{rem:super2}.
\end{remark}

\section{Fluid analysis and asymptotic optimality}
\label{sec:fluid}
In this section, we present a fluid formulation of the restless bandit
problem and show that its optimal fluid cost provides a lower bound on
the cost in the original stochastic model. Based on the optimal fluid
solution, we then derive a set of priority policies for the original
stochastic model that we prove to be asymptotically optimal.

In Section~\ref{sec:fluidmodel}, we introduce the fluid control
problem. In Section~\ref{sec:priority}, we define the set of priority
policies and the asymptotic optimality results can be found in
Section~\ref{sec:asympt}.

\subsection{Fluid control problem and lower bound}
\label{sec:fluidmodel}




The fluid control problem arises from the original stochastic model by
taking into account only the mean drifts. For a given control $u(t)$,
let $x_{j,k}^{u,a}(t)$ denote the amount of class-$k$ fluid in
state $j$ under action $a$ at time $t$ and let
$x_{j,k}^{u}(t)=x_{j,k}^{u,0}(t)+x_{j,k}^{u,1}(t)$ be the amount of
class-$k$ fluid in state $j$. The dynamics is then given by
%
\begin{eqnarray}
\label{eq:gen} \frac{\mathrm{d}x_{j,k}^{u}(t)}{\mathrm{d}t} &=& \lambda_k p_k(j) +
\sum_{a=0}^{1} \sum
_{i=1, i\neq j }^{J_k} x^{u,a}_{i,k}(t)
q_k(j|i,a)\nonumber\\
&&{} - \sum_{a=0}^{1}
\sum_{i=0, i\neq j }^{J_k} x^{u,a}_{j,k}(t)
q_k(i|j,a)
\\
&=& \lambda_k p_k(j) +\sum
_{a=0}^{1} \sum_{i=1 }^{J_k}
x^{u,a}_{i,k}(t) q_k(j|i,a),\nonumber
\end{eqnarray}
where the last step follows from $q_k(j|j,a):=-\sum_{i=0, i\neq j}^{J_k}q_k(i|j,a)$.
The constraint on the total amount of active fluid is given by
\[
\sum_{k=1}^K\sum
_{j=1}^{J_k} x^{u,1}_{j,k}(t)\leq
\alpha\qquad \mbox{for all } t\geq0.
\]
We are interested in finding an optimal equilibrium point of the fluid
dynamics that minimizes the holding cost. Hence, we pose the following
linear optimization problem:
%
\begin{eqnarray}
\mathrm{(LP)}\quad &&\min_{(x_{j,k}^a)} \sum_{k=1}^K
\sum_{j=1}^{J_k}\sum
_{a=0}^{1} C_k(j,a)
x^a_{j,k}
\nonumber
\\
&&\mbox{s.t.}\qquad 0 = \lambda_k p_k(j) +\sum
_{a=0}^{1} \sum_{i=1 }^{J_k}
x^a_{i,k} q_k(j|i,a)\qquad \forall j,k,
\label{eq:dif0}
\\
&&\hphantom{\mbox{s.t.}\qquad} \sum_{k=1}^K\sum
_{j=1}^{J_k} x^1_{j,k} \leq
\alpha,\label
{eq:Alpha}
\\
&& \hphantom{\mbox{s.t.}\qquad}\sum_{j=1}^{J_k} \sum
_{a=0}^{1}x_{j,k}^a =
x_k(0)\qquad \mbox{if } \lambda_k =0, \forall k,
\label{eq:27}
\\
&& \hphantom{\mbox{s.t.}\qquad}x^a_{j,k}\geq0\qquad \forall j,k,a, \label{eq:ja2}
\end{eqnarray}
where the constraint \eqref{eq:27} can be seen as follows: if $\lambda
_k=0$, then $q_k(0|i,a)=0$ for all $i$. Hence, from \eqref{eq:gen} we
obtain $\sum_{j=1}^{J_k} \frac{\mathrm{d}}{\mathrm{d}t} x^u_{j,k}(t)=0$.

We denote by $x^*$ an optimal solution of the above problem (LP),
assuming it exists. For a fixed population, an optimal solution depends
on $x_k(0)$. However, for ease of notation, this dependency is not
stated explicitly. We further denote the optimal value of the (LP) by
\[
v^*\bigl(x(0)\bigr): =\sum_{k=1}^K\sum
_{j=1}^{J_k}\sum
_{a=0}^1 C_{k}(j,a)
x^{*,a}_{j,k}.
\]
%

We can now state some results concerning the optimization problem (LP).
The proof of the first lemma may be found in Appendix \ref{appa}.

\begin{lemma}
\label{15151}
If there exists a policy $\pi$ such that the process $X^\pi(t)$ has a
unique invariant probability distribution with finite first moments,
then the feasible set of \textup{(LP)} is nonempty and $v^*(x)<\infty$, for
any $x$.
\end{lemma}

As a consequence of Lemma~\ref{15151}, we get a necessary
condition under which there exists a policy that makes the system
stable and has finite first moments.

\begin{corollary}
If there exists a policy $\overline\pi$ such that the system is stable
with finite first moments, then
\[
\sum_{k=1}^K \sum
_{j=1}^{J_k} y_{j,k}^{*1} \leq
\alpha,
\]
with
$y^{*} :=\arg\min\{\sum_{k=1}^K\sum_{j=1}^{J_k} x^1_{j,k} : x \mbox{
satisfies } \eqref{eq:dif0}, \eqref{eq:27} \mbox{ and } \eqref{eq:ja2}
\}$.
\end{corollary}

\begin{pf}
Assume there exists a policy $\overline\pi$
such that the process $X^{\overline\pi}(t)$ has a unique invariant
probability distribution with finite first moments. By Lemma~\ref
{15151}, the feasible set of (LP) is nonempty. That is, there
exists an $(x^a_{j,k})$ such that \eqref{eq:dif0}, \eqref{eq:27} and
\eqref{eq:ja2} hold and $\sum_{k=1}^K\sum_{j=1}^{J_k} x^1_{j,k} \leq
\alpha$. Hence, by definition of the optimal solution $y^*$ we obtain
$\sum_{k=1}^K \sum_{j=1}^{J_k} y_{j,k}^{*1} \leq\sum_{k=1}^K\sum_{j=1}^{J_k} x^1_{j,k} \leq\alpha$. This completes the proof.
\end{pf}

The optimal solution of the fluid control problem (LP) serves as a
lower bound on the cost of the original stochastic optimization
problem, see the following lemma. The proof can be found in Appendix \ref{appb}.

\begin{lemma}
\label{25252}
%
For a fixed population of bandits, we have that for any policy $\pi$,
%
\begin{equation}
\label{eq:liminf_merged} 
V_-^\pi(x) \geq
v^*(x).
\end{equation}

For a dynamic population of bandits, relation \eqref{eq:liminf_merged}
holds if
\begin{itemize}
%
\item policy $\pi$ is stable, or,
\item policy $\pi$ is (mean) rate-stable and $C_{k}(j,a)>0$, for all $j,k,a$.
\end{itemize}
%
\end{lemma}


\subsection{Priority policies}
\label{sec:priority}

A priority policy is defined as follows. There is a predefined priority
ordering on the states each bandit can be in. At any moment in time, a
priority policy makes active a maximum number of bandits being in the
states having the highest priority among all the bandits present. In
addition, the policy can prescribe that certain states are never made active.

We now define a set of priority policies $\Pi^*$ that will play a key
role in the paper.
The priority policies are derived from (the) optimal equilibrium
point(s) $x^*$ of the (LP) problem: for a given equilibrium
point $x^*$, we consider all priority orderings such that the states
that in equilibrium are never passive ($x_{j,k}^{*,0}=0$) are of higher
priority than states that receive some passive action ($x^{*,0}_{j, k}>0$).
In addition, states that in equilibrium are both active and passive
($x^{*,0}_{j,k} \cdot x^{*,1}_{j,k}>0$) receive higher priority than
states that are never active ($x^{*,1}_{j,k}=0$).
Further, if the full capacity is not used in equilibrium (i.e.,
$\sum_k \sum_j x^{*,1}_{j,k}<\alpha$), then the states that are never active
in equilibrium are never activated in the priority ordering.
The set of priority policies $\Pi^*$ is formalized in the definition below.

\begin{definition}[(Set of priority policies)]
\label{def:Pix}
We define
\[
X^*:= \bigl\{x^*: x^* \mbox{ is an optimal solution of (LP) with
$x_k(0)=X_k(0)$}\bigr\}.
\]
The set of priority policies $\Pi^*$ is defined as
\[
\Pi^*:=\bigcup_{x^*\in X^*} \Pi\bigl(x^*\bigr),
\]
%
where $\Pi(x^*)$ is the set of all priority policies that satisfy the
following rules:
\begin{longlist}[1.]
\item[1.] A class-$k$ bandit in state $j$ with $x^{*,1}_{j,k}>0$ and
$x^{*,0}_{j,k}=0$ is given higher priority than a class-$\tilde k$
bandit in state $\tilde j$ with $x^{*,0}_{\tilde j,\tilde k}>0$.
\item[2.] A class-$k$ bandit in state $j$ with $x^{*,0}_{j, k}>0$ and
$x^{*,1}_{ j, k} >0$ is given higher priority than a class-$\tilde k$
bandit in state $\tilde j$ with $x^{*,0}_{\tilde j,\tilde k}>0$ and
$x^{*,1}_{\tilde j,\tilde k} =0$.
\item[3.] If $\sum_{k=1}^K \sum_{j=1}^{J_k} x^{*,1}_{j,k}<\alpha$, then any
class-$k$ bandit in state $j$ with $x^{*,1}_{j,k}=0$ and
$x^{*,0}_{j,k}>0$ will \textit{never} be made active.
\end{longlist}
\end{definition}

We emphasize that in order to define the set of priority policies $\Pi
^*$, we do not require the bandits to be indexable, as defined in
Definition~\ref{def:indexability}. This is in contrast to the
definition of Whittle's index policy, which is only well defined in the
case the system is indexable.
We note that Whittle's index policy is included in $\Pi^*$ for
indexable systems as will be proved in Section~\ref{sec:asW}.

If there exists a policy such that the system is stable and has finite
first moments, then the feasible set of (LP) is nonempty (Lemma~\ref
{15151}), and hence the set $\Pi^*$ is nonempty.
Note that the set $\Pi^*$ can consist of more than one policy.
When selecting a policy it might be of practical importance to aim for
a policy that is \emph{robust} in the arrival characteristics, the
number of bandits that can be made active and the number of bandits in
each class.

\begin{definition}[(Robust policy)]
\label{def:robust}
A priority policy is called \emph{robust} if the priority ordering does
not depend on $\alpha$, $\lambda_k$ and $X_k(0)$, $k=1,\ldots,K$.
\end{definition}

In the case the system is indexable, Whittle's index policy is a robust
element of $\Pi^*$; see Section~\ref{sec:Wheur}.
For a nonindexable system, the set $\Pi^*$ might no longer contain a
robust policy. In Section~\ref{sec:nonindex}, we explain how to select
in that case priority policies from the set $\Pi^*$.


Before continuing, we first give an example of Definition~\ref{def:Pix}.

\begin{example}
Assume $K=2$ and $J_k=2$. Let $x^*$ be such that for class 1 we have
$x_{1,1}^{*,0}=0$, $x_{2,1}^{*,0}=4$, $x_{1,1}^{*,1}=3$,
$x_{2,1}^{*,1}=1$ and for class 2 we have $x_{1,2}^{*,0}= 2$,
$x_{2,2}^{*,0}= 0$, $x_{1,2}^{*,1}= 0$, $x_{2,2}^{*,1}= 5$ and $\alpha=10$.
The priority policies associated to $x^*$ in the set $\Pi(x^*)$, as
defined in Definition~\ref{def:Pix}, satisfy the following rules:
By point 1: class-1 bandits in state 1 and class-2 bandits in state 2
are given the highest priority.
By point 3: since $x_{1,1}^{*,1}+ x_{2,1}^{*,1} + x_{1,2}^{*,1}+
x_{2,2}^{*,1} =9<\alpha$, class-2 bandits in state 1 are never made active.
Let the pair $(j,k)$ denote a class-$k$ bandit in state $j$. The set
$\Pi(x^*)$ contains two policies: either give priority according to
$(1,1)\succ(2,2) \succ(2,1)$ or give priority according to
$(2,2)\succ(1,1) \succ(2,1)$. In neither policy, state $(1,2)$ is
never made active.
\end{example}

\begin{remark}[(Multi actions)]
\label{rem:super1}
In this remark, we explain how to define the set of priority policies
$\Pi^*$ in the case of multiple actions per bandit.
Similar to the nonrestless bandit problem (see Remark~\ref
{rem:multi}), we are interested in priority policies such that if a
class-$k$ bandit in state $j$ is chosen to be active, it will always be
made active in a fixed mode $a_{k}(j)\in\{0, 1,2, \ldots, A_{k}(j)\}$.
We therefore need to restrict the set $X^*$ to optimal solutions of
(LP) that satisfy $x_{j,k}^{*,a }x_{j,k}^{*,\tilde a}=0$, for all $a,
\tilde a \in\{ 1, \ldots, A_{k}(j)\}$. The latter condition implies
that for all activation modes $a=1,\ldots, A_{k}(j)$ one has
$x_{j.k}^{*,a}=0$, with the exception\vspace*{1pt} of at most one active mode,
denoted by $a_{k}(j)$.
The set $\Pi(x^*)$ is then defined as in Definition~\ref{def:Pix},
replacing the action $a=1$ by $a=a_{k}(j)$.
All results obtained in Section~\ref{sec:fluid} remain valid [replacing
$a=1$ by $a=a_{k}(j)$].
\end{remark}

\begin{remark}
\label{rem:NM}
In \cite{BNM00}, a heuristic is proposed for the multi-class restless
bandit problem for a fixed population of bandits: the so-called
primal--dual heuristic. This is defined based on the optimal (primal and
dual) solution of an LP problem corresponding to the discounted-cost
criterion. In fact, if the primal--dual heuristic would have been
defined based on the problem (LP), it can be checked that it satisfies
the properties of Definition~\ref{def:Pix}, and hence is included in
the set of priority policies $\Pi^*$.
\end{remark}


In order to prove asymptotic optimality of a policy $\pi^*\in\Pi^*$,
as will be done in Section~\ref{sec:asympt}, we investigate its fluid dynamics.
Denote by $S^{\pi^*}_{k}(j)$ the set of pairs $(i,l), i=1,\ldots, J_l,
l=1,\ldots, K$, such that class-$l$ bandits in state $i$ have higher
priority than class-$k$ bandits in state $j$ under policy $\pi^*$.
Denote by $I^{\pi^*}$ the set of all states that will never be made
active under policy $\pi^*$.
The fluid dynamics under policy $\pi^*$ can now be written as follows:
\begin{eqnarray}
\label{eq:ode}  \frac{\mathrm{d} x^{\pi^*}_{j,k}(t)}{\mathrm{d}t}& =& \lambda_k p_k(j)
+ \sum_{a=0}^1 \sum
_{i=1 }^{J_k} x^{\pi^*,a}_{i,k}(t)
q_k(j|i,a),\nonumber
\\
\eqntext{\mbox{with } x^{\pi^*,1}_{j,k}(t)=\min \biggl( \biggl(\alpha-
\displaystyle\sum_{(i,l) \in S^{\pi^*}_{k}(j)} x^{\pi^*}_{i,l}(t)
\biggr)^+, x^{\pi
^*}_{j,k}(t) \biggr)\qquad\mbox{if $(j,k)\notin
I^{\pi^*}$,}}
\\[-8pt]
\\[-8pt]
\nonumber
x^{\pi^*,1}_{j,k}(t)&=&0\quad\mbox{ if $(j,k)\in
I^{\pi^*}$,}
\nonumber
\\
 x^{\pi^*,0}_{j,k}(t)&=&x^{\pi^*}_{j,k}(t)-x_{j,k}^{\pi^*,1}(t).
\nonumber
\end{eqnarray}
It follows directly that an optimal solution $x^*$ of (LP) is an
equilibrium point of the process $x^{\pi^*}(t)$.

\begin{lemma}
Let $\pi^*\in\Pi^*$ and let $x^*$ be a point such that $\pi^*\in
\Pi
(x^*)$. Then $x^*$ is an equilibrium point of the process $x^{\pi
^*}(t)$ as defined in \eqref{eq:ode}.
\end{lemma}

\begin{pf}
Since $x^*$ is an optimal solution of (LP), it follows directly from
the definition of $\Pi(x^*)$ that $x^*$ is an equilibrium point of the
process $x^{\pi^*}(t)$.
%
\end{pf}

In order to prove asymptotic optimality of a policy $\pi^*$, we will
need that the equilibrium point $x^*$ is in fact a global attractor of
the process $x^{\pi^*}(t)$, that is, all trajectories converge to
$x^*$. This is not true in general, which is why we state it as a
condition for a policy to satisfy. In Section~\ref{sec:global}, we will
further comment on this condition.


\begin{condition}
\label{cond:better}
Given an equilibrium point $x^*\in X^*$ and a policy $\pi^*\in\Pi
(x^*)\subset\Pi^*$. The point $x^*$ is a global attractor of the
process $x^{\pi^*}(t)$. That is, for any initial point, 
the process $x^{\pi^*}(t)$ converges to $x^*$.
\end{condition}

\subsection{Asymptotic optimality of priority policies}
\label{sec:asympt}

In this section, we present the asymptotic optimality results for the
set of priority policies $\Pi^*$.
In particular, we obtain that the priority policies minimize the \emph
{fluid-scaled} average holding cost.

We will consider the restless bandit problem in the following
fluid-scaling regime: we scale by $r$ both the arrival rates and the
number of bandits that can be made active. That is, class-$k$ bandits
arrive at rate $\lambda_k \cdot r$, $k=1,\ldots, K$, and $\alpha
\cdot
r$ bandits can be made active at any moment in time.
We let $X^{r}_{j,k}(0)= x_{j,k}\cdot r$, with $x_{j,k}\geq0$.
For a given policy $\pi$,\vspace*{1pt} we denote by $X^{r,\pi,a}_{j,k}(t)$ the
number of class-$k$ bandits in state $j$ experiencing action $a$ at
time $t$ under scaling parameter $r$.

We make the important observation that the set of policies $\Pi^*$ is
invariant to the scaling parameter. This follows since an optimal
solution of (LP) scales with the parameter $r$: if $x^*$ is an optimal
solution, then so is $x^*r$ for the (LP) with parameters $\alpha\cdot
r$, $x(0)\cdot r$ and $\lambda_k \cdot r$. By Definition~\ref{def:Pix},
the set of priority policies does therefore not depend on $r$.

We will be interested in the process after the fluid scaling, that is,
space is scaled linearly with the parameter $r$, $\frac{X^{r,\pi
,a}_{j,k}(t)}{r}$.
%
We further define for a given initial state $x$,
\[
V_-^{r,\pi}(x) := \liminf_{T\to\infty}\frac{1}{T}
\ME_{r\cdot
x} \Biggl(\int_0^T \sum
_{k=1}^K \sum
_{j=1}^{J_k}\sum_{a=0}^1
C_{k}(j,a)\frac
{X_{j,k}^{r,\pi,a}(t)}{r} \,\mathrm{d}t \Biggr)
\]
and
\[
V_+^{r,\pi}(x) := \limsup_{T\to\infty}\frac{1}{T}
\ME_{r\cdot
x} \Biggl(\int_0^T \sum
_{k=1}^K \sum
_{j=1}^{J_k}\sum_{a=0}^1
C_{k}(j,a) \frac
{X_{j,k}^{r,\pi,a}(t)}{r} \,\mathrm{d}t \Biggr).
\]
If $V_-^{r,\pi}(x)=V_+^{r,\pi}(x)$ for all $x$, then we define
$V^{r,\pi
}(x):= V_+^{r,\pi}(x)$.

Our goal is to find policies that minimize the cost of the stochastic
model after fluid scaling.
We therefore call a policy $\pi^*$ asymptotically optimal when the
fluid-scaled version of (\ref{eq:avopt}) holds.

\begin{definition}[(Asymptotic optimality)]
A policy $\pi^*$ is asymptotically optimal if
\[
\limsup_{r\to\infty} V_+^{r,\pi^*}(x) \leq\liminf
_{r\to\infty} V_-^{r,\pi}(x) \qquad\mbox{for all $x$ and all
policies $\pi\in G$},
\]
where $G$ is a set of admissible policies.
\end{definition}

In our asymptotic optimality result, the set $G$ will consist of all
policies for the fixed population of bandits, while it will consists of
all policies that are stable, rate-stable or mean rate-stable for the
dynamic population of bandits; see Proposition~\ref{45454}.



In order to prove asymptotic optimality of priority policies in the set
$\Pi^*$, we need the following technical condition.

\begin{condition}
\label{cond:twee}
Given a policy $\pi^*\in\Pi^*$.
\begin{longlist}[(a)]
%
\item[(a)] The process $\frac{X^{r,\pi^*}(t)}{r}$ has a unique invariant
probability distribution $p^{r,\pi^*}$, which has a finite first
moment, for all $r$.
\item[(b)] The family $\{p^{r,\pi^*},r\}$ is tight.
\item[(c)] The family $\{p^{r,\pi^*},r\}$ is uniform integrable.
\end{longlist}
\end{condition}

For a fixed population of bandits, the state space of $X^{r,\pi^*}(t)$
is finite, hence conditions (b) and (c) are satisfied. A sufficient
condition for Condition \ref{cond:twee}(a) to hold is the Markov
process $X^{r,\pi^*}(t)$ to be unichain, for any $r$, \cite{Tijms2003}.

For a dynamic population of bandits, we present a large class of
restless bandit problems for which Condition \ref{cond:twee} is satisfied.
More precisely, we consider problems in which bandits that are kept
passive will eventually leave the system.
For many real-life situations, this assumption arises naturally. For
example, customers that become impatient and abandon the queue/system,
companies that go bankrupt, perishable items, etc.
The proof of the proposition may be found in Appendix \ref{appc}.

\begin{proposition}
\label{35353}
Assume 
that the state $0$ is positive recurrent for a class-$k$ bandit that is
kept passive.
For any priority policy $\pi$ for which $X^{r,\pi}(t)$ is irreducible,
Condition \ref{cond:twee} is satisfied.
%
\end{proposition}

Another class of problems satisfying Condition \ref{cond:twee} would be
those in which only active bandits are allowed in the system, that is,
$q_k(0|i,0)=\infty$, for all $k,i$. This could describe for example the
hiring process where new candidates are modeled by new arriving
bandits, room occupation in a casualty departments where patients
require direct attention, or a loss network.
When $q_k(0|i,0)=\infty$, for all $k, i$, at most $\alpha$ bandits are
present in the system, hence due to the finite state space,
Condition \ref{cond:twee} follows directly from a unichain assumption.

We can now state the asymptotic optimality result.



\begin{proposition}
\label{45454}
For a given policy $\pi^*\in\Pi(x^*)\subset\Pi^*$, assume
Conditions \ref{cond:better} and~\ref{cond:twee} are
satisfied. Then
\[
\lim_{r\to\infty} V^{r,\pi^*}(x)=v^*(x)\qquad\mbox{for any $x$.}
\]
%
In particular, we have
\[
\liminf_{r\to\infty} V_-^{r,\pi}(x)\geq\lim
_{r\to\infty} V^{r,\pi
^*}(x)\qquad \mbox{for any $x$ and any policy $
\pi\in G$,}
\]
%
where for the fixed population of bandits $G$ consists of all policies,
and for the dynamic population of bandits
\begin{itemize}
%
\item$G$ is the set of all stable policies $\pi$, or,
\item$C_{k}(j,a)>0$, for all $j,k,a$ and $G$ is the set of all
rate-stable and mean rate-stable policies.
\end{itemize}
%
\end{proposition}

The proof may be found in Appendix \ref{appd} and consists of the following
steps: Given a policy $\pi^*\in\Pi(x^*)$, we show that the
fluid-scaled steady-state queue length vector converges to $x^*$. Since
$x^*$ is an optimal solution of the fluid control problem (LP) with
$x(0)=x$ and has cost value $v^*(x)$, this implies that the
fluid-scaled cost under policy $\pi^*$ converges to $v^*(x)$. Since
$v^*(x)$ serves as a lower bound on the average cost, this allows us to
conclude for asymptotic optimality of the priority policy $\pi^*$.

\section{Whittle's index policy}
\label{sec:whittle}

In Section~\ref{sec:asympt}, we showed that priority policies inside
the set $\Pi^*$ are asymptotically optimal. In this section, we will
derive that Whittle's index policy is included in this set of policies
$\Pi^*$.

In Section~\ref{sec:rel}, we first define Whittle's index policy. In
Sections \ref{sec:asWNoArrivals} and~\ref{sec:asW}, we then give
sufficient conditions under which Whittle's index policy is
asymptotically optimal, both in the case of a fixed population of
bandits, and in the case of a dynamic population of bandits, respectively.

\subsection{Relaxed-constraint optimization problem and Whittle's indices}
\label{sec:rel}

Whittle's index policy was proposed by Whittle \cite{Whittle88} as an
efficient heuristic for the multi-class restless bandit problem.
Each bandit is assigned a Whittle's index, which is a function of the
state the bandit is in. Whittle's index policy activates those bandits
having currently the highest indices.
In this section, we will describe how these Whittle's indices are derived.

In order to define Whittle's indices, we consider the following
optimization problem:
Find a stationary and Markovian policy that minimizes
%
\begin{equation}
\label{eq:ma1} \mathbb{C}_x^{f} \Biggl( \sum
_{k=1}^K \sum_{j=1}^{J_k}
\sum_{a=0}^{1} C_{k}(j,a)
X_{j,k}^{\pi,a}(\cdot) \Biggr) \qquad\mbox{with } f\in\{\mathrm{av}, \beta
\},
\end{equation}
under the constraint \eqref{eq:M},
where
%
\begin{equation}
\label{eq:costav} \mathbb{C}_x^{\mathrm{av}}\bigl( Y(\cdot) \bigr):=
\limsup_{T\to\infty} \frac
{1}{T}\mathbb {E}_x
\biggl(\int_0^T Y(t)\,\mathrm{d}t \biggr),
\end{equation}
represents the average-cost criterion and
\[
\mathbb{C}_x^{\beta}\bigl(Y(\cdot) \bigr):=
\mathbb{E}_x \biggl(\int_0^\infty
\mathrm{e}^{-\beta t} Y(t)\,\mathrm{d}t \biggr),
\]
$\beta>0$, represents the discounted-cost criterion.
The objective as stated in Section~\ref{sec:model} is the average-cost
criterion. In Section~\ref{sec:asW}, it will become clear why we need
to introduce here the discounted-cost criterion as well.

\subsubsection{Relaxed-constraint optimization problem}
\label{subsec:relax}
The restless property of the bandits makes the above described
optimization problem often infeasible to solve. Instead, Whittle \cite
{Whittle88} proposed to study the so-called \emph{relaxed-constraint
optimization problem}, which is defined as follows:
find a policy that minimizes \eqref{eq:ma1}
under the relaxed constraint
%
\begin{equation}
\label{eq:rel_2} \mathbb{C}_x^{f} \Biggl( \sum
_{k=1}^K \sum_{j=1}^{J_k}
X_{j,k}^{\pi
,1}(\cdot) \Biggr) \leq\alpha(f),
\end{equation}
with $\alpha(\mathrm{av})= \alpha$ and $\alpha(\beta)= \int_0^\infty\alpha
\mathrm{e}^{-\beta t}\,\mathrm{d}t= \alpha/\beta$ for $\beta>0$.
That is,
the constraint that at most $\alpha$ bandits can be made active at any
moment in time is replaced by its time-average or discounted
version, \eqref{eq:rel_2}.
Hence, the cost under the optimal policy of the relaxed-constraint
optimization problem provides a lower bound on the cost for any policy
that satisfies the original constraint.

In standard restless bandit problems, the constraint \eqref{eq:rel_2}
needs to be satisfied in the strict sense, that is, with an ``$=$'' sign.
In this paper, we allow however strictly less than $\alpha$ bandits to
be active at a time. In order to define Whittle's indices,
we therefore introduce so-called \emph{dummy bandits}. That is, besides
the initial population of bandits, we assume there are $\alpha(f)$
additional bandits that will never change state. We denote the state
these bandits are in by $B$ and the cost of having a dummy bandit in
state $B$ is $C_B(a)=0$, $a=0,1$. The introduction of these $\alpha(f)$
dummy bandits allows to reformulate the relaxed-constraint problem as
follows: minimize \eqref{eq:ma1} under the relaxed constraint
%
\begin{equation}
\label{eq:rel_3} \mathbb{C}_x^{f} \bigl(X_{B}^{\pi,1}(
\cdot) \bigr) + \mathbb {C}_x^{f} \Biggl( \sum
_{k=1}^K \sum_{j=1}^{J_k}
X_{j,k}^{\pi,1}(\cdot) \Biggr) = \alpha(f).
\end{equation}
This constraint is equivalent to \eqref{eq:rel_2} since, for a given
set of active bandits, activating additional dummy bandits does not
modify the behavior of the system.

Using the Lagrangian approach, we write the relaxed-constraint problem
[minimize \eqref{eq:ma1} under constraint \eqref{eq:rel_3}] as the
problem of finding a policy $\pi$ that minimizes
%
\begin{eqnarray}
\label{eq:ma9} &&\sum_{k=1}^K \sum
_{j=1}^{J_k} \mathbb{C}_x^{f}
\bigl( C_{k}(j,0) X_{j,k}^{\pi,0}(\cdot)
+C_{k}(j,1) X_{j,k}^{\pi,1}(\cdot) +\nu
X_{j,k}^{\pi,1}(\cdot) \bigr)
\nonumber
\\[-8pt]
\\[-8pt]
\nonumber
&&\qquad{}+ \mathbb{C}_x^{f}
\bigl( \nu X_{B}^{\pi
,1}(\cdot) \bigr).
\end{eqnarray}
The Lagrange multiplier $\nu$ can be viewed as the cost to be paid per
active bandit. From Lagrangian relaxation theory, we have that there
exists a value of the Lagrange multiplier $\nu$ such that the
constraint (\ref{eq:rel_3}) is satisfied.

Since there is no longer a common constraint for the bandits,
problem \eqref{eq:ma9} can be decomposed into several \emph
{subproblems}, one for each bandit: for each class-$k$ bandit the
subproblem is to minimize
%
\begin{equation}
\label{eq:subproblem} \mathbb{C}^{f} \bigl( C_k
\bigl(J_k(\cdot),A^\pi_k(\cdot)\bigr) +\nu
\mathbf{1}_{( A^\pi_k(\cdot)=1)} \bigr),
\end{equation}
where $J_k(t)$ denotes the state of a class-$k$ bandit at time $t$ and
$A^\pi_k(t)$ denotes the action chosen for the class-$k$ bandit under
policy $\pi$.
We take as convention that $J_k(t)=0$ and $A_k(t)=0$ if the bandit is
not present (or no longer present) in the system at time $t$ and set
$C_k(0,0)=0$.
For each dummy bandit, the problem is to minimize
%
\begin{equation}
\label{eq:prr} \nu \mathbb{C}^{f} ( \mathbf{1}_{( A^\pi_B(\cdot)=1)} ),
\end{equation}
with $A^\pi_B(t)$ the action chosen for the dummy bandit at time $t$
under policy $\pi$.

We can now define Whittle's index.

\begin{definition}[(Whittle's index)]
For a given optimization criterion $f$, we define Whittle's index $\nu
_k^f(j)$ as the least value of $\nu$ for which it is optimal in \eqref
{eq:subproblem} to make the class-$k$ bandit in state $j$ passive.

Similarly, we define the index $\nu_{B}^f$ as the least value of $\nu$
for which it is optimal in \eqref{eq:prr} to make a dummy bandit passive.
\end{definition}

Indexability is the property that allows to characterize an optimal
policy for the relaxed optimization problem.

\begin{definition}[(Indexability)]
\label{def:indexability}
A bandit is \emph{indexable} if the set of states in which passive is
an optimal action in \eqref{eq:subproblem}, denoted by $D(\nu)$,
increases in $\nu$. That is,
$\nu'<\nu$ implies $D(\nu')\subset D(\nu)$.
\end{definition}

%
We note that the dynamics of a bandit in state $B$ is independent of
the action chosen. Since $\nu$ represents the cost to be paid when
active, it will be optimal in \eqref{eq:prr} to make a bandit in
state $B$ passive if and only if $\nu\geq0$. As a consequence, a~dummy
bandit is always indexable and $\nu_{B}^f =0$.

We call the problem indexable if all bandits are indexable.
Note that whether or not a problem is indexable can depend on the
choice for $f$ (and $\beta$). We refer to \cite{NM07} for a survey on
indexability results.
In particular, \cite{NM07} presents sufficient conditions for a
restless bandit to be indexable and provides a method to calculate
Whittle's indices. Sufficient conditions for indexability can also be
found in \cite{LZ10,W07}.

If the bandit problem is indexable, an optimal policy for the
subproblem \eqref{eq:subproblem} is then such that the class-$k$ bandit
in state $j$ is made active if $\nu_{k}^f(j)> \nu$, is made passive if
$\nu_{k}^f(j) < \nu$, and any action is optimal if $\nu_{k}^f(j)=\nu
$, \cite{Whittle88}.

An optimal solution to \eqref{eq:ma1} under the relaxed
constraint \eqref{eq:rel_3} is obtained by setting $\nu$ at the
appropriate level $\nu^*$ such that \eqref{eq:rel_3} is satisfied. A
class-$k$ bandit in state $j$ is then made active if $\nu_{k}^f(j)>\nu
^*$, and kept passive if $\nu_{k}^f(j)<\nu^*$. When a class-$k$ bandit
is in a state $j$ such that $\nu_{k}^f(j)=\nu^*$, one needs to
appropriately randomize the action in this state such that the relaxed
constraint \eqref{eq:rel_3} is satisfied, \cite{WW90,Whittle88}.
In the case $\nu^*=0$, we take the convention that the randomization is
done among the bandits in state $B$ [possible since there are exactly
$\alpha(f)$ dummy bandits], while any class-$k$ bandit in a state $j$
with $\nu_{k}^f(j)=0$ is kept passive.

Since $\nu_B^f=0$, a dummy bandit has higher priority than a class-$k$
bandit in state $j$ with $\nu_{k}^f(j)\leq0$.
Together with constraint \eqref{eq:rel_3} and the fact that there are
$\alpha(f)$ dummy bandits, we conclude that any class-$k$ bandit in
state $j$ with $\nu_{k}^f(j)\leq0$ is kept passive in the relaxed
optimization problem. In particular, this implies that $\nu^*\geq0$.

\subsubsection{Whittle's index policy as heuristic}
\label{sec:Wheur}
The optimal control for the relaxed problem is not feasible for the
original optimization problem having as constraint that at most $\alpha
$ bandits can be made active \emph{at any moment in time}.
Whittle \cite
{Whittle88} therefore proposed the following heuristic:

\begin{definition}[(Whittle's index policy)]
For a given optimization criterion~$f$, Whittle's index policy
activates the $\alpha$ bandits having currently the highest \emph
{nonnegative} Whittle's index value $v_k^f(j)$.
In case different states have the same value for the Whittle index, an
arbitrary fixed priority rule is used.
We denote Whittle's index policy by $\nu^{f}$.
\end{definition}

If $\nu_{k}^f(j)<\nu_{l}^f(i)$, then a class-$l$ bandit in state $i$ is
given higher priority than a class-$k$ bandit in state $j$ under
Whittle's index policy.
Analogously to the optimal solution of the relaxed optimization
problem, a class-$k$ bandit in state $j$ with $\nu^f_{k}(j) \leq0$
will never be made active under Whittle's index policy. It can
therefore happen that strictly less than $\alpha$ bandits are made
active, even though there are more than $\alpha$ bandits present.

Whittle's indices result from solving \eqref{eq:subproblem}.
Since the latter does not depend on $\alpha$, $\lambda_k$, and
$X_k(0)$, we can conclude that Whittle's index policy is a \emph
{robust} policy; see Definition~\ref{def:robust}.
In the next two sections, we will prove that Whittle's index policy is
asymptotically optimal, both for the static and dynamic population.

\begin{remark}[(Multi-actions)]
\label{rem:super2}
In this remark, we define Whittle's index policy in the case of
multiple actions. For that we need to assume a stronger form of
indexability: There is an index $v_{k}^f(j)$ and an activation mode
$a_{k}^f(j)$ such that an optimal solution of \eqref{eq:subproblem} is
to make a class-$k$ bandit in state $j$ active in mode $a_{k}^f(j)$ if
$\nu<\nu_{k}^f(j)$ and to keep it passive if $\nu>\nu_{k}^f(j)$.
Whittle's index rule is then defined as in Section~\ref{subsec:relax},
replacing the action $a=1$ by $a=a_{k}^f(j)$.

If the restless bandit problem satisfies this stronger form of
indexability, then one can reduce the multi-action problem to the
single-action problem and hence all asymptotic optimality results as
obtained in Sections \ref{sec:asWNoArrivals} and~\ref{sec:asW}
will be valid [replacing action $a=1$ by $a=a_{k}(j)$].
\end{remark}

\subsection{Asymptotic optimality for a fixed population of bandits}
\label{sec:asWNoArrivals}

In this section, we consider a fixed population of indexable bandits and
show that Whittle's index policy, defined for the time-average cost
criterion $f=\mathrm{av}$, is asymptotically optimal.


We will need the following assumption, which was also made in \cite{WW90}.

\begin{assumption}
\label{as:unichain}
For every $k$, the process describing the state of a class-$k$ bandit
is unichain,
regardless of the policy employed.
\end{assumption}

The next proposition shows that Whittle's index policy is included in
the set of priority policies $\Pi^*$.
The proof can be found in Appendix \ref{appe}.

\begin{proposition}
\label{56565}
Consider a fixed population of bandits. If Assumption~\ref{as:unichain}
holds and if the restless bandit problem is indexable for the
average-cost criterion, then there is an $x^*\in X^*$ such that
Whittle's index policy $\nu^{\mathrm{av}}$ is included in the set $\Pi
(x^*)\subset\Pi^*$.
\end{proposition}

We can now conclude that Whittle's index policy is asymptotically optimal.

\begin{corollary}
\label{cor:fixed}
Consider a fixed population of bandits. If the assumptions of
Proposition~\ref{56565} are satisfied and if
Condition \ref{cond:better} holds for Whittle's index policy $\nu
^{\mathrm{av}}$, then
\[
\lim_{r\to\infty} V^{r, \nu^{\mathrm{av}}}(x) \leq\liminf
_{r\to\infty}V_-^{r,\pi}(x),
\]
%
for any $x$ and any policy $\pi$.
\end{corollary}

\begin{pf}
From Propositions \ref{45454} and~\ref
{56565}, we obtain the desired result.
\end{pf}


 The above corollary was previously proved by Weber and Weiss
in \cite{WW90} for the case of symmetric bandits, that is, $K=1$. We note
that the assumptions made in~\cite{WW90} in order to prove the
asymptotic optimality result are the same as the ones in Corollary~\ref
{cor:fixed}.

The proof technique used in Weber and Weiss \cite{WW90} is different
from the one used here. In \cite{WW90}, the cost under an optimal
policy is lower bounded by the optimal cost in the relaxed problem and
upper bounded by the cost under Whittle's index policy.
By showing that both bounds converge to the same value, the fluid
approximation, the asymptotic optimality of Whittle's index policy is concluded.
Obtaining a lower bound for a dynamic population does not seem
straightforward. This is why we undertook in this paper a different
proof approach that applies as well for a dynamic population; see
Section~\ref{sec:asW}.

\subsection{Asymptotic optimality for a dynamic population of bandits}
\label{sec:asW}

In this section, we will introduce an index policy for the dynamic
population of bandits, based on Whittle's indices, and show it to be
asymptotically optimal.
More precisely, we show the index policy to be included in the set of
asymptotically optimal policies $\Pi^*$, as obtained in Section~\ref{sec:asympt}.

Recall that our objective is to find a policy that asymptotically
minimizes the average-cost criterion \eqref{eq:costav}.
We do however not make use of Whittle's index policy $\nu^{\mathrm{av}}$ for the
following reason: Consider a class-$k$ bandit and the relaxed
optimization problem \eqref{eq:subproblem}, with $f=\mathrm{av}$. Any policy
that makes sure that the class-$k$ bandit leaves after a finite amount
of time has an average cost equal to zero and is hence an optimal solution.
In order to derive a nontrivial index rule, the authors of \cite
{AEJ10,AJN10} consider instead the Whittle indices corresponding to
the discounted-cost criterion ($f=\beta$, $\beta>0$).
An index rule for the average-cost criterion is then obtained by
considering the limiting values as $\beta\downarrow0$.
We propose here the same.
%
%
For a given class $k$, let $\beta_l\downarrow0$ be some subsequence
such that the limit
\[
\nu_{k}^{\mathrm{lim}}(j):=\lim_{l\to\infty}
\nu_{ k}^{\beta_l}(j)
\]
exists, for all $j=1,\ldots,J_k$. The limit can possibly be equal to
$\infty$.
The index policy $\nu^{\mathrm{lim}}$ activates the $\alpha$ bandits having
currently the highest \emph{nonnegative} index value $\nu_k^{\mathrm{lim}}(j)$.
In this section, we will show asymptotic optimality of $\nu^{\mathrm{lim}}$. In
order to do so, we will need that class-$k$ bandits are indexable under
the $\beta_l$-discounted cost criterion, for $l$ large enough. In
addition, we will need the following assumption on the model parameters.

\begin{assumption}
\label{as:bounded2}
For all $k=1,\ldots, K$, the set of optimal solutions of the linear program
\begin{eqnarray*}
&&\min_x \sum_{j=1}^{J_k}
\bigl(C^0_{j,k} x^0_{j,k} +
C^1_{j,k} x^1_{j,k} +\nu
x^1_{j,k} \bigr)
\\
&&\mbox{s.t.}\qquad 0= \lambda_k p_k({0,j}) + \sum
_{a=0}^1 \sum_{i=1}^{J_k}
x^a_{i,k} q_k(j|i,a)\qquad \forall j,
\\
&&\hphantom{\mbox{s.t.}\qquad} x^a_{j,k}\geq0\qquad \forall j, a,
\end{eqnarray*}
is bounded when $\nu>0$.
\end{assumption}

We note that this assumption is always satisfied if $C_{k}(j,0)>0$ and
$C_{k}(j,1)\geq0$, for all $j,k$, since $x_{j,k}^{*,1}$ and
$x^{*,0}_{j,k}$ are upper bounded by the cost value of a feasible
solution divided by $\nu+C_{k}(j,1)>0$ and $C_{k}(j,0)>0$, respectively.



The proposition below shows that Whittle's index policy $\nu^{\mathrm{lim}}$ is
included in the set of priority policies $\Pi^*$. The proof can be
found in Appendix \ref{appe}.

\begin{proposition}
\label{75757}
Consider a dynamic population of bandits.
For a given class $k$, let $\beta_l\downarrow0$ be some subsequence
such that the limit
\[
\nu_{k}^{\mathrm{lim}}(j):=\lim_{l\to\infty}
\nu_{ k}^{\beta_l}(j)
\]
exists, for all $j=1,\ldots,J_k$.
If Assumption~\ref{as:bounded2} holds and if the discounted restless
bandit problem is indexable for $\beta_l \leq\overline\beta$, with
$0<\overline\beta< 1$,
then there is an $x^*\in X^*$ such that Whittle's index policy $\nu
^{\mathrm{lim}}$ is included in the set $\Pi(x^*)\subset\Pi^*$.
\end{proposition}

We can now conclude for asymptotic optimality of Whittle's index policy
$\nu^{\mathrm{lim}}$.

\begin{corollary}
\label{cor:dynamic}
Consider a dynamic population of bandits.
If the assumptions of Proposition~\ref{75757} are satisfied and
if Conditions \ref{cond:better} and~\ref{cond:twee} hold for
Whittle's index policy $\nu^{\mathrm{lim}}$, then
%
\begin{eqnarray}
&&\lim_{r\to\infty\label{eqLopt_rate2}} V^{r,\nu^{\mathrm{lim}}}(x) \leq\liminf
_{r\to\infty} V_-^{r,\pi}(x) \qquad\mbox{for all $x$ and any policy
$\pi\in G$},
\end{eqnarray}
%
where
\begin{itemize}
\item$G$ consists of all stable policies, or
\item$C_{k}(j,a)>0$, for all $j,k,a$ and $G$ consists of all
rate-stable and mean rate-stable policies.
\end{itemize}
\end{corollary}

\begin{pf}
The result follows directly from Propositions \ref{45454} and~\ref{75757}.
\end{pf}



The above result for the dynamic population shows that the heuristic
$\nu^{\mathrm{lim}}$, which is based on a model without arrivals, is in fact
nearly optimal in the presence of arrivals.
In addition, Whittle's index policy $\nu^{\mathrm{lim}}$ is robust, that is, it
does not depend on the arrival characteristics of new bandits or on the
exact number of bandits that can be made active.

\begin{remark}[(Multi-actions)] In order to define $\nu^{\mathrm{lim}}$ in the
case of multiple actions per bandit, we need to assume that, for $\beta
_l$ small enough, the stronger form of indexability (defined in
Remark~\ref{rem:super2}) holds. In addition, the optimal activation
mode for a class-$k$ bandit in state $j$, denoted by $a_{k}^{\beta
_l}(j)$, cannot depend on $\beta_l$, that is, $a_{k}^{\beta_l}(j)=a_{k}(j)$.
\end{remark}

\section{On the global attractor property}
\label{sec:examples}
\label{sec:global}

In Proposition~\ref{45454}, asymptotic optimality of priority
policies in the set $\Pi^*$ was proved under
the global attractor property (Condition \ref{cond:better}) and a
technical condition (Condition \ref{cond:twee}).
In this section, we further discuss the global attractor property.
The latter is concerned with the process $x^*(t)$, defined by the
ODE (\ref{eq:ode}), to have a global attractor.
We recall that in~\cite{WW90} the same global attractor property was
required in order to prove asymptotic optimality of Whittle's index
policy for a fixed population of symmetric bandits ($K=1$). In
addition, the authors of \cite{WW90} presented an example for which
Whittle's index policy is not asymptotically optimal (and hence, does
not satisfy the global attractor property).

For a fixed population of symmetric indexable bandits, the global
attractor property was proved to always hold under Whittle's index
policy if a bandit can be in at most three states ($J=3$); see \cite{WW91}.
However, in general no sufficient conditions are available in order for
$x^*$ to be a global attractor of $x^{\pi^*}(t)$.
A necessary condition was provided in \cite{WW90}, Lemma~2, where for a
fixed population of symmetric bandits it was proved that indexability
is necessary in order for Whittle's index policy to satisfy the global
attractor property, for any value of $\alpha$ and $x(0)$.
We emphasize that when the system is nonindexable, there can still
exist priority policies in $\Pi^*$ (possibly nonrobust) that satisfy
the global attractor property.

The asymptotic optimality result of Whittle's index policy for the case
$K=1$, \cite{WW90}, has been cited extensively. The global attractor
property is often verified only numerically.
Note that in the context of mean field interaction models, convergence
of the stationary measure also relies on a global attractor assumption
of the corresponding ODE; see, for example, \cite{BL08}.
In a recent paper, the authors of \cite{OES12} proved asymptotic
optimality of Whittle's index policy for a very specific model with
only two classes of bandits (fixed population of bandits) under a
recurrence condition. The latter condition replaced the global
attractor condition, however, the authors needed as well to resort to
numerical experiments in order to verify this recurrence condition.

In the remainder of this section, we describe the necessity of the
global attractor property and the technical challenges in the case this
condition is not satisfied.

Optimal fluid control problems have been widely studied in the
literature in order to obtain asymptotically optimal policies for the
stochastic model. In the context of this paper, the fluid control
problem related to our results would be to find the optimal
control $u^*(t)$ that minimizes
%
\begin{equation}
\label{eq:optd} \lim_{T\to\infty}\frac{1}{T}\int
_0^T \sum_{k=1}^K
\sum_{j=1}^{J_k} \sum
_{a=0}^1 C_{k}(j,a)
x_{j,k}^{u,a}(t) \,\mathrm{d}t,
\end{equation}
where the dynamics of $x_{j,k}^{u,a}(t)$ is described by \eqref{eq:gen}.
The optimal control $u^*(t)$ is then to be translated back to the
stochastic model in such a way that it is asymptotically optimal.
When stating the global attractor property, the above is exactly what
we have in mind. In fact, instead of solving this transient fluid
control problem, we directly consider an optimal equilibrium point of
the fluid model and propose a priority policy based on this equilibrium
point. When the global attractor property is satisfied, this implies
that the optimal equilibrium point is indeed reached by the associated
strict priority control, and hence this priority control solves \eqref
{eq:optd}.

When for any $\pi^*\in\bigcup_{x^*\in X^*} \Pi(x^*)=\Pi^*$ the global
attractor property is not satisfied, this means that there does not
exist a priority control $u(t)=\pi^*\in\Pi(x^*)$ such that the fluid
process $x^{\pi^*}(t)$ converges to $x^*$. In that case, we can be in
either one of the following two situations: (1) There exists a control
$u^*(t)$ for which the process $x^{u^*}(t)$ does have as global
attractor $x^*\in X^*$, where $X^*$ was defined as the set of optimal
equilibrium points. This control $u^*(t)$ might not be of priority
type. (2) There does not exist any control that has a global attractor
$x^*\in X^*$.
In the latter case, the optimal control $u^*(t)$ can be such that the
process $x^{u^*}(t)$ behaves cyclically or shows chaotic behavior, or
the process converges to a nonoptimal equilibrium point.
Hence, in the case Condition \ref{cond:better} is not satisfied, in
both situations (1) and (2), one needs to determine the exact
transient behaviour of the optimal control of \eqref{eq:optd},
$u^*(t)$, which in its turn needs to be translated back to the
stochastic model.
We leave this as subject for future research.

\section{Case study: A multi-server queue with abandonments}
\label{sec:abopt}

In this section, we study a multi-class multi-server system with
impatient customers, the multi-class $\mathit{M/M/S+M}$ system. This is an
example of a restless bandit problem with a dynamic population. We will
derive a robust priority policy that is in the set $\Pi^*$ and show
that it satisfies the two conditions needed in order to conclude for
asymptotic optimality.

The impact of abandonments has attracted considerable interest from the
research community, with a surge in recent years. To illustrate the
latter, we can mention the recent Special Issue on abandonments in
Queueing Systems \cite{HP13} and the survey paper \cite{DH12} on
abandonments in a many-server setting


We consider a multi-class system with $S$ servers working in parallel.
At any moment in time, each server can serve at most one customer.
Class-$k$ customers arrive according to a Poisson process with rate
$\lambda_k>0$ and require an exponentially distributed service with
mean $1/\mu_k<\infty$. Server $s$, $s=1,\ldots, S$ works at speed 1.
Customers waiting (being served) abandon the queue after an
exponentially distributed amount of time with mean $1/\theta_k$
($1/\tilde\theta_k$), with $\theta_k>0$, $\tilde\theta_k \geq0$, for
all $k$. Having one class-$k$ customers waiting in the queue (in
service) costs $c_k$ ($\tilde c_k$) per unit of time. Each abandonment
of a waiting class-$k$ customer (class-$k$ customer being served) costs
$d_k$ ($\tilde d_k$).
We are interested in finding a policy $\pi$ that minimizes the long-run
average cost
\[
\limsup_{T\to\infty} \frac{1}{T}\sum
_{k=1}^K \ME_x \biggl(\int
_0^T \bigl( c_k
X_k^{\pi,0}(t) + \tilde c_k X_k^{\pi,1}(t)
\bigr) \,\mathrm{d}t + d_k R_k^\pi(T) +\tilde
d_k \tilde R_k^\pi(T) \biggr),
\]
where $X_k^{\pi,0}(t)$ [$X_k^{\pi,1}(t)$] denotes the number of
class-$k$ customers in the queue (in service) at time $t$ and $R_k^\pi
(t)$ [$\tilde R_k^\pi(t)$] denotes the number of abandonments of
waiting class-$k$ customers (class-$k$ customers being served) in the
interval $[0,t]$.

Representing each customer in the queue (in service) by a passive
(active) bandit, the problem can be addressed within the framework of a
restless bandit model with the following parameters:\vspace*{1pt} $J_k=1$,
$q_k(0|1,0)=\theta_k>0$, $q_k(0|1,1)= \mu_k +\tilde\theta_k$, $C_k(1,0) =
c_k+d_k \theta_k$, $C_k(1,1) = \tilde c_k + \tilde d_k
\tilde\theta_k$,\vspace*{1pt}
$k=1,\ldots, K$, and $\alpha=S$, where we used that $\ME_x(R_k^\pi(T))
=\theta_k\ME_x(\int_0^T X_k^{\pi,0}(t)\,\mathrm{d}t)$ and\break $\ME
_x(\tilde
R_k^\pi(T)) =\tilde\theta_k\ME_x(\int_0^T X_k^{\pi,1}(t)\,\mathrm{d}t)$.
A bandit can only be in two states (state 0 or state 1), hence
indexability follows directly (for any choice of $\beta$).

We now define an index policy that we will prove to be included in the
set $\Pi^*$.
For each class $k$, we set
\[
\iota_k := q_k(0|1,1) \biggl( \frac{C_k(1,0) }{q_k(0|1,0)} -
\frac
{C_k(1,1)}{q_k(0|1,1)} \biggr). 
\]
The index policy $\iota$ is then defined as follows: At any moment in
time serve (at most) $S$ customers present in the system that have the
highest, strictly positive, index values, $\iota_k$. If a customer
belongs to a class that has a negative index value, then this customer
will never be served.

Before continuing, we first give an interpretation of the index $\iota_k$.
The term $1/q_k(0|1,a)$ is the time it takes until a bandit under
action $a$ leaves the system. Hence, $C_k(1,a) /q_k(0|1,a)$ is the cost
for applying action $a$ on a class-$k$ bandit until it leaves the
system. The difference $ \frac{C_k(1,0) }{q_k(0|1,0)} - \frac
{C_k(1,1)}{q_k(0|1,1)}$ is the reduction in cost when making a
class-$k$ bandit active (instead of keeping him passive), so that the
index $\iota_k$ represents the reduction in cost per time unit when
class $k$ is made active.
Also note that the index rule $\iota$ does not depend on the arrival
rate of the customers or the number of servers present in the system,
hence it is a robust rule; see Definition~\ref{def:robust}.

By solving the LP problem corresponding to the multi-server queue with
abandonments, we obtain in Proposition~\ref{prop:ab} that the index
policy $\iota$ is included in $\Pi^*$.
%
\begin{proposition}
\label{prop:ab}
Policy $\iota$ is contained in the set $\Pi^*$.

In addition, when $\iota_1 > \iota_2 > \cdots> \iota_K$, policy
$\iota
$ coincides with Whittle's index policy $\nu^{\mathrm{lim}}$.
\end{proposition}

\begin{pf}
For the multi-class multi-server system with abandonments, the linear
program (LP) is given by
%
\begin{eqnarray}\label{eq:bal}
&&\min_x \sum_k
\bigl(c_k x^0_k +\tilde c_k
x^1_k + d_k \theta_k
x^0_k + \tilde d_k \tilde
\theta_k x^1_k\bigr),
\nonumber
\\
&&\mbox{s.t.}\qquad 0 = \lambda_k - \mu_k
x^1_k - \theta_k x^0_k
-\tilde\theta_k x^1_k,
\\
&& \hphantom{\mbox{s.t.}\qquad}\sum_{k=1}^K x_k^1
\leq S\quad \mbox{and}\quad x_k^0, x_k^1
\geq 0.
\nonumber
\end{eqnarray}
Equation \eqref{eq:bal} implies $x_k^0 = \frac{\lambda_k - (\mu_k
+\tilde\theta_k)x_k^1}{\theta_k}$. Hence, the above linear program is
equivalent to solving
\begin{eqnarray*}
&&\max_x \sum_k \biggl(
(c_k +d_k\theta_k )\frac{\mu_k +\tilde
\theta
_k}{\theta_k} -
\tilde c_k - \tilde d_k \tilde\theta_k
\biggr)x_k^1,
\\
&&\mbox{s.t.}\qquad \sum_{k=1}^K
x_k^1 \leq S\quad \mbox{and}\quad 0\leq x_k^1
\leq\frac
{\lambda_k}{\mu_k+\tilde\theta_k}.
\end{eqnarray*}
The optimal solution is to assign maximum values to those $x_k^{1}$
having the highest values for $\iota_k=(c_k +d_k\theta_k )\frac{\mu_k
+\tilde\theta_k}{\theta_k} - \tilde c_k - \tilde d_k \tilde\theta_k$,
with $\iota_k>0$, until the constraint $\sum_k x_k^1\leq S$ is
saturated. Denote this optimal solution by $x^*$.
Assume the classes are ordered such that $\iota_1 \geq\iota_2 \geq
\cdots\geq\iota_K$.\vspace*{1pt}
Hence, one can find an $l$ such that: (1) for all $k<l$ it holds that
$x_k^{*,1}=\frac{\lambda_k}{\mu_k +\tilde\theta_k}$, and hence
$x_k^{*,0}=0$, (2) for $k=l$ it holds that $0\leq x_l^{*,1}\leq\frac
{\lambda_l}{\mu_l +\tilde\theta_l}$, and hence $x_l^{*,0}\geq0$,
(3) and for all $k>l$ it holds that $x_k^{*,1}=0$. This gives that the
index policy $\iota$ is included in the set $\Pi(x^*) \subset\Pi^*$;
see Definition~\ref{def:Pix}.

When $\iota_1 > \iota_2 > \cdots> \iota_K$, it follows directly that
$\iota$ is the unique policy that is in the set $\Pi^*$ for any value
of $S$ or $\lambda_k$. Since Whittle's index policy is by definition
robust and is in the set $\Pi^*$ (Proposition~\ref{75757}), we
obtain that $\iota$ and Whittle's index policy have the same priority ordering.
\end{pf}

Note that the $\mathit{M/M/S+M}$ system belongs to the class of problems as
described in Proposition~\ref{35353}. Hence, Condition \ref
{cond:twee} is satisfied. The global attractor property follows
from \cite{AGS10}, where this property was proved for a slightly
different model. We therefore have the following optimality result for
the index policy $\iota$.

\begin{proposition}
\label{cor:ab}
Consider a system with $Sr$ servers working in parallel and arrival
rates $\lambda_k r$, $k=1,\ldots, K$. The index policy $\iota$ is
asymptotically optimal as $r\to\infty$, that is, for any $x$ and any
policy $\pi$,
\begin{eqnarray*}
&&\lim_{r\to\infty} \lim_{T\to\infty}\frac{1}{T}
\ME_x \Biggl( \int_0^T \sum
_{k=1}^K \biggl((c_k+d_k
\theta_k) \frac
{X^{r,\iota,0}_{k}(t)}{r} + (\tilde c_k+\tilde
d_k \tilde\theta_k) \frac{X^{r,\iota,1}_{k}(t)}{r} \biggr)
\,\mathrm{d}t \Biggr)
\\
&&\qquad \leq\liminf_{r\to\infty} \liminf_{T\to\infty}
\frac{1}{T} \ME_x \Biggl( \int_0^T
\sum_{k=1}^K \biggl((c_k+d_k
\theta_k) \frac
{X^{r,\pi,0}_{k}(t)}{r }\\
&&\qquad\quad{}+ (\tilde c_k+\tilde
d_k \tilde\theta_k) \frac
{X^{r,\pi,1}_{k}(t)}{r} \biggr)
\,\mathrm{d}t \Biggr).
\end{eqnarray*}
\end{proposition}

\begin{pf}
In Proposition~\ref{prop:ab}, we showed that $\iota$ is included in
$\Pi
(x^*)$, with $x^*$ as given in the proof of Proposition~\ref{prop:ab}.
In Appendix \ref{apph}, we prove that the process $x^{\iota}(t)$, as defined
in \eqref{eq:ode}, has the point $x^*$ as a unique global attractor,
that is, Condition \ref{cond:better} is satisfied.
From Proposition~\ref{35353}, we obtain that Condition \ref
{cond:twee} is satisfied.
Further, note that any policy $\pi$ gives a stable system, since
$\theta
_k>0$ for all $k$.
Together with Proposition~\ref{45454}, we then obtain that
the index policy $\iota$ is asymptotically optimal.
 \end{pf}

%
%

\begin{remark}[(Existing results in literature)]
In \cite{LAV14}, a single-server queue with abandonments has been
studied. Whittle's index policy was there derived by modeling the
system as a fixed population of restless bandits: each bandit
representing a class and the state of a bandit representing the number
of customers in the queue. The latter implies that $J_k=\infty$, for
all $k$, hence it does not fall inside the framework of this paper. The
results obtained in \cite{LAV14} apply to general holding cost
functions. In the case of linear holding costs, as considered in this
section, the index rule as derived in \cite{LAV14} coincides with
policy $\iota$.
We further note that even though in \cite{LAV14} the arrival
characteristics are taken into account when calculating Whittle's
indices, the final result is independent on the arrival
characteristics. For nonlinear holding cost, this is no longer the case.

For the case $\tilde c_k=0, \tilde\theta_k=0$ and $c_k + d_k \theta
_k>0$, the asymptotic optimality of the policy $\iota$ in a
multi-server setting has previously been proved in \cite{AGS10,AGS11}.
Note that in this setting the performance criterion is the weighted
number of customers present in the queue. If $\sum\lambda_k/\mu_k>S$,
that is, the overload situation, the fluid-scaled cost $v^*(S)$ will be
nonzero, and hence the optimality result is useful. This is not the
case when $\sum\lambda_k/\mu_k< S$, the underload setting, as was also
observed in \cite{AGS10,AGS11}: in underload we have for any
nonidling policy $x_k^{*,0}=0$, $\forall k$, see equation \eqref
{eq:a1}, which together with $\tilde c_k=0$ implies $v^*(S)=0$, that
is, in equilibrium the cost is zero for any nonidling policy.
In \cite{LAV13}, the \emph{transient behavior} of the fluid model has
been studied for the underload setting. It was shown that the optimal
transient fluid control is in fact a state-dependent strategy and hence
no longer a strict priority policy.

For a discrete-time model with one server and $\tilde\theta_k=0,
\tilde c_k=c_k>0$, Whittle's index $\nu_{k}^{\mathrm{lim}}$ has been derived
in \cite{AJN10}. This index $\nu_{k}^{\mathrm{lim}}$ coincides with the
Whittle's index $\iota_k$ for the continuous-time model.
In this setting, the fluid-scaled cost is always strictly positive:
$v^*(S)=0$ would imply that $x^*_k=0$, however, this contradicts with
equation \eqref{eq:bal}, which would read $0=\lambda_k$. Hence, the
asymptotic optimality result applies to both the underload and overload regime.
\end{remark}

\section{Nonindexable restless bandits}
\label{sec:nonindex}

The set of priority policies, $\Pi^*$, consists of more than one
policy, and hence, it is not direct which priority policy to choose.
For an indexable restless bandit problem, Whittle's index policy is
inside the set $\Pi^*$ and is robust, that is, it does not depend on
$\alpha, \lambda_k, X_k(0)$, $k=1,\ldots, K$. This is therefore an
obvious choice, and Whittle's index policy has been extensively tested
numerically for different applications and shown to perform well; see,
for example, \cite{ALJZK09,AGNK03,AJN10,AEJV13,EM04,GKO09,GM02,LAV14,LZ10,N07,RBCK08} and the examples in the book \cite{GGW11}.
In this section, we therefore focus our attention on nonindexable
restless bandits. In Section~\ref{sec:nonindex_serie}, we describe how
to select a priority policy from the (possibly large) set of priority
policies $\Pi^*$ and in Section~\ref{sec:nonindex_example} their
performance is numerically evaluated outside the asymptotic regime.

\subsection{Policy selection}
\label{sec:nonindex_serie}

In this section, we describe how to select priority policies from the
set $\Pi^*$. In order to do so, we will need the following technical
lemma that gives a characterization of an optimal solution of the (LP)
problem. Note that this lemma is valid for both indexable and
nonindexable examples.
We refer to Appendix \ref{appg} for the proof.

\begin{lemma}
\label{888885}
In the case of a dynamic population of bandits, assume that the set of
optimal solutions of \textup{(LP)} is bounded and either $p_k(j)>0$, for all
$j,k$, or $C_k(j,0) >0$, for all $j,k$.

For either a fixed or dynamic population of bandits, there exists at
least one optimal solution of \textup{(LP)}, $x^*$, such that
$x_{j,k}^{*,0}x_{j,k}^{*,1}>0$ for at most one pair $(j,k)$.
\end{lemma}

The assumption that the set of optimal solutions of (LP) is bounded is
always satisfied if $C_k(j,0) > 0$, for all $j,k$. This follows since
$x_{j,k}^{*,1}\leq\alpha$ and $x^{*,0}_{j,k}\leq ( \overline C-
\sum_{j,k}C_k(j,1) x_{j,k}^{*,1} )/C_k(j,0) <\infty$, with
$\overline C<\infty$ the cost value of a feasible solution.

In the remainder of this section, we will write $\Pi^*(\alpha)$ instead
of $\Pi^*$ to emphasize the dependence on $\alpha$, that is, the number
of bandits that can be simultaneously made active. In the case of
indexable bandits, there exists priority policies that are inside $\Pi
^*(\alpha)$, for all $\alpha$, for example, Whittle's index policy. In
general, this is not the case for nonindexable bandits.
Below we therefore describe how one can select priority policies from
the set $\Pi^*(\alpha)$ as $\alpha$ changes.

From Lemma~\ref{888885}, we have that, for a fixed $\alpha$, there
exists at least one optimal solution of (LP), $x^*(\alpha)$, such that
$x_{j,k}^{*,0}(\alpha)x_{j,k}^{*,1}(\alpha)>0$, for at most one pair $(j,k)$.
In particular, we can define $0=\alpha_0< \alpha_1<\alpha_2<\cdots
<\alpha_M$ and $\alpha_{M+1}=\infty$, such that for a given interval
$[\alpha_i,\alpha_{i+1})$ the binding constraints of the (LP) and the
basis of an optimal solution do not change. Hence, there are pairs
$(j_i,k_i)$ and sets $H_i, L_i$ and $\tilde L_i$ such that, for any
$\alpha\in[\alpha_i,\alpha_{i+1})$, it holds that
\begin{eqnarray*}
x_{j,k}^{*,0}(\alpha) &= &0 \quad\mbox{and}\quad
x_{j,k}^{*,1}(\alpha) \geq0\qquad \mbox{for all } (j,k)\in
H_i,
\\
x_{j_i,k_i}^{*,0}(\alpha)&\geq&0\quad \mbox{and}\quad
x_{j_i,k_i}^{*,1}(\alpha ) \geq0,
\\
x_{j,k}^{*,0}(\alpha) &\geq&0 \quad\mbox{and}\quad
x_{j,k}^{*,1}(\alpha) = 0\qquad \mbox{for all } (j,k)\in
L_i,
\\
x_{j,k}^{*,0}(\alpha) &= &0 \quad\mbox{and}\quad
x_{j,k}^{*,1}(\alpha) = 0\qquad \mbox{for all } (j,k)\in\tilde
L_i,
\end{eqnarray*}
and either $\sum_{k=1}^K \sum_{j=1}^{J_k} x_{j,k}^{*,1}(\alpha)=
\alpha
$ or $\sum_{k=1}^K \sum_{j=1}^{J_k} x_{j,k}^{*,1}(\alpha)< \alpha$.

When choosing a priority policy from the set $\Pi^*(\alpha)$, we
propose to choose the same policy for any $\alpha\in[\alpha_i,
\alpha
_{i+1})$. This policy is chosen in the following way:
\begin{itemize}
\item
Class-$k$ bandits in state $j$ with $(j,k)\in H_i$ receive highest priority.
\item Class-$k_i$ bandits in state $j_i$ receive lower priority than
class-$k$ bandits in state $j$ with $(j,k)\in H_i$.
\item For class-$k$ bandits in state $j$ with $(j,k)\in L_i$, we have
to distinguish between two situations:
\begin{itemize}[(ii)]
\item[(i)] if $\sum_{k=1}^K \sum_{j=1}^{J_k} x_{j,k}^{*,1}(\alpha)<
\alpha$, that is, there is capacity left unused, then any class-$k$
bandit in state $j$, with $(j,k)\in L_i$, will never be made active.
\item[(ii)] if\vspace*{1pt} $\sum_{k=1}^K \sum_{j=1}^{J_k} x_{j,k}^{*,1}(\alpha)=
\alpha$, then the capacity constraint is binding. We will allow bandits
in the set $L_i$ to be made active only if this would have happened
when there would have been more capacity $\alpha$ available.
Hence, a~class-$k$ bandit in state $j$, $(j,k)\in L_i$, receives lower
priority than a class-$\tilde k$ bandit in state $\tilde j$, $(\tilde
j, \tilde k)\in H_i\cap\{j_i,k_i\}$, if there is an $n>i$ such that
$(j,k)\in H_n \cap\{j_n,k_n\}$. If there does not exist such $n$, then
such bandits are never made active.
\end{itemize}
\item Class-$k$ bandits in state $j$ with $(j,k)\in\tilde L_i$ are
never made active.
\end{itemize}
It is left open how to set the priority ordering within the high
priority states $H_i$ and the low priority states $L_i$. One way would
be to chose the priorities such that the priority ordering changes
minimally as $\alpha$ changes to other intervals.


\subsection{Performance evaluation}
\label{sec:nonindex_example}

We now turn our attention to a particular nonindexable example and
numerically evaluate the selection method as explained in the previous section.
We took the continuous-time version of the example given in~\cite{NM07},
Section~2.2. We consider a fixed population of bandits, and
each bandit can be in three states. The cost structure is given by
\[
\bigl(C(1,0), C(2,0), C(3,0)\bigr)=(-0.458, -0.5308, -0.6873)
\]
and
\[
\bigl(C(1,1), C(2,1), C(3,1)\bigr)=(-0.9631,-0.7963, -0.1057).
\]
The transition matrices $Q^0=(q(j|i,0))_{i,j}$ and
$Q^1=(q(j|i,1))_{i,j}$ are given by
%
\begin{eqnarray}
Q^0&=& %
\pmatrix{ -0.8098& 0.4156 & 0.3942 \vspace*{2pt}
\cr
0.5676& -0.5809 &0.0133 \vspace*{2pt}
\cr
0.0191 &0.1097& -0.1288 } %
\quad\mbox{and}
\nonumber
\\[-8pt]
\\[-8pt]
\nonumber
 Q^1&=& %
\pmatrix{ -0.2204 & 0.0903 &0.1301
\vspace*{2pt}
\cr
0.1903 &-0.8137& 0.6234 \vspace*{2pt}
\cr
0.2901& 0.3901& -0.6802
} %
.
\end{eqnarray}

Our aim in this section is to numerically evaluate the performance of
priority policies in $\Pi^*$ outside the asymptotic regime. In
particular, we evaluate the performance when $\alpha=1$, that is, at
most one bandit can be made active at a time, and we let the number of
bandits, $X(0)$, vary.

For a given value of $\alpha$ and $X(0)$, the set $\Pi^*$ consists of
more than one policy. Before presenting the numerical results, we will
therefore first describe the priority policies we considered using the
selection method as given in the previous section.
In Table~\ref{table_sol}, one can find the structure of an optimal
basic solution of (LP), obtained numerically, when fixing $\alpha=1$
and letting the number of bandits present in the system, $x(0)$,
increase. We note that equivalently we could have taken $x(0)=\bar x$
fixed and let $\alpha$ decrease, simply by a change of variable in the
(LP) problem.

\begin{table}
\tabcolsep=0pt
\caption{Optimal basic solutions for the \textup{(LP)} problem}
\label{table_sol}
\begin{tabular*}{\textwidth}{@{\extracolsep{\fill}}lcccc@{}}
\hline
\multicolumn{1}{@{}l}{$\bolds{\alpha=1}$} &
\multicolumn{4}{c@{}}{\textbf{Optimal basic solution}}
\\
\hline
$ 1\leq x(0)\leq2.4 $ & $x_1^{*,0}=0, x_1^{*,1}>0$ & $x_2^{*,0}>0,
x_2^{*,1}=0$ & $x_3^{*,0}>0, x_3^{*,1}=0$& \eqref{eq:Alpha} not
binding\\
$2.4\leq x(0)\leq3.6$ & $x_1^{*,0}=0, x_1^{*,1}>0$ & $x_2^{*,0}>0,
x_2^{*,1}>0$ & $x_3^{*,0}>0, x_3^{*,1}=0$& \eqref{eq:Alpha} binding\\
$ 3.6\leq x(0)\leq7.36 $ & $x_1^{*,0}>0, x_1^{*,1}>0$ &
$x_2^{*,0}=0, x_2^{*,1}>0$ & $x_3^{*,0}>0, x_3^{*,1}=0$& \eqref
{eq:Alpha} binding\\
$7.36\leq x(0)$ & $x_1^{*,0}>0, x_1^{*,1}=0$ & $x_2^{*,0}=0,
x_2^{*,1}>0$ & $x_3^{*,0}>0, x_3^{*,1}=0$& \eqref{eq:Alpha} binding\\
\hline
\end{tabular*}
\end{table}

\begin{table}[b]
\caption{Selected priority policies}
\label{table_pol}
\begin{tabular*}{\textwidth}{@{\extracolsep{\fill}}lccc@{}}
\hline
& \multicolumn{1}{c}{\textbf{Priority ordering}} &
\multicolumn{1}{c}{\textbf{Always passive}} &
\multicolumn{1}{c@{}}{\textbf{Name of policy}}
\\
\hline
$ x(0)= 1,2 $ & 1 & 2, 3 & prio1 \\
$ x(0)= 3 $ & $1 \succ2$& 3 & prio12 \\
$ x(0)\geq4$ & $2 \succ1$& 3 & prio21 \\
\hline
\end{tabular*}
\end{table}

We can now characterize the priority policies; see also Table~\ref{table_pol}.
Consider $x(0)=1$ or $x(0)=2$. In that case, we derive from Table~\ref{table_sol} that a bandit in state 1 receives priority (by
Definition~\ref{def:Pix}). Since the constraint \eqref{eq:Alpha} is not
binding, a bandit in state 2 or 3 will never be made active (by
Definition~\ref{def:Pix}). Hence, $\Pi(x^*)$ consists of the policy
that only makes bandits active in state 1. This policy is referred to
as ``prio1''.
Now consider $x(0)=3$. Then Definition~\ref{def:Pix} prescribes that,
for any policy in $\Pi(x^*)$, state 1 has strict priority over state 2,
and state 3 is either never made active, or has lowest priority. Note
that for smaller values of $x(0)$ (equivalent to considering higher
values of $\alpha$), state 3 is not made active either. Hence, as
explained in the previous section we choose to keep bandits in state 3
passive, that is, we focus on the policy ``prio12''.
Now consider $4\leq x(0)\leq7$. Then Definition~\ref{def:Pix}
prescribes that, for any policy in $\Pi(x^*)$, state 2 has strict
priority over state 1, and state 3 is either never made active, or has
lowest priority. Note that state 3 is never made active for $x(0)<4$.
Hence, as explained in the previous section, we chose to do the same
for $4\leq x(0)\leq7$, that is, we focus on policy ``prio21''.
Now consider $x(0)\geq8$. Then Definition~\ref{def:Pix} prescribes
that, for any policy in $\Pi(x^*)$, state 2 has strict priority and
that states 1 and 3 are either never made active or have lowest
priority. For smaller values of $x(0)$, class 1 is made active, while
class 3 is never made active. Hence, as explained in the previous
section, we chose to do the same for $x(0)\geq8$, that is, we focus on
policy ``prio21''.

Any policy gives a unichain Markov chain, hence Condition \ref
{cond:twee} is satisfied. We therefore have that any priority policy in
$\Pi^*$ that satisfies the global attractor property, as in
Condition \ref{cond:better}, is asymptotically optimal.
Numerically, we evaluated the global attractor property and found the following:
for $x(0)=1$, the policy prio1 has $x^*$ as global attractor, while
policies prio12 and prio123, which also belong to $\Pi^*$, converge to
a nonoptimal equilibrium point.
For $x(0)=3$, there does not exist a priority policy that converges to
$x^*$. For example, the fluid dynamics under prio12 converges to an
equilibrium where state 1 is sometimes passive (and state 2 and 3 are
never active), while the optimal point $x^*$ never makes state 1 passive.
For $4\leq x(0)\leq7$, the set $\Pi^*$ consists of the policies prio21
and prio213, both of them have $x^*$ as global attractor.
For $x(0)\geq8$, the set $\Pi^*$ consists of prio2, prio21 and
prio213, all of them have $x^*$ as global attractor.

\begin{figure}

\includegraphics{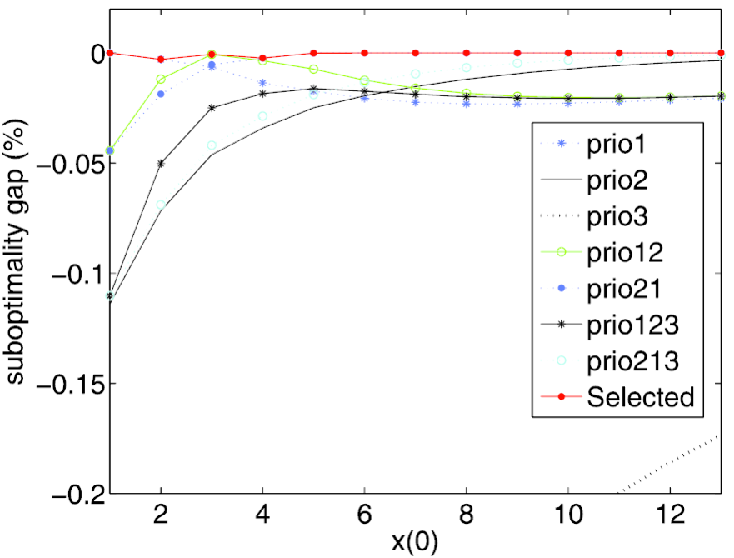}

\caption{Sub-optimality gap of priority policies for nonindexable example.}
\label{fig:NM_nonindexleft}
\end{figure}

We have numerically evaluated the performance of the priority policies
as described in Table~\ref{table_pol} against both the optimal policy
(obtained numerically by value iteration) and against other priority
policies. In Figure~\ref{fig:NM_nonindexleft}, we plot the relative
sub-optimality gap (in \%) for the different policies when $\alpha=1$
and let the number of bandits, $X(0)=x(0)$, vary on the horizontal
axis. The line referred to as ``selected'' plots for each given $x(0)$
the selected priority policy as given in Table~\ref{table_pol}. We
observe that these selected policies always have the smallest
sub-optimality gap.

Prio123 and prio213 are inside the class of asymptotically optimal
policies, $\Pi^*$, for $X(0)=3$ and $X(0)\geq4$, respectively,
however, the selection process, as described in Section~\ref{sec:nonindex_serie}, does not select these policies. In fact, we
observe that prio123 and prio213 are outperformed by our selected
priority policies. Below we will see that this sub-optimality gap can
be made arbitrarily large.

\begin{figure}

\includegraphics{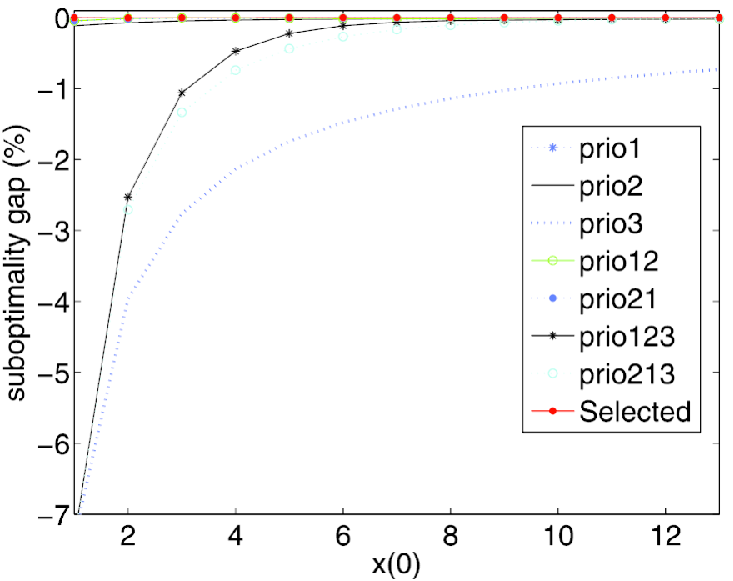}

\caption{Sub-optimality gap of priority policies for nonindexable example.}
\label{fig:NM_nonindexright}
\end{figure}

The difference in performance between different priority policies is
not that large in this example. For other instances of nonindexable
bandits, including dynamic populations, the differences can be larger though.
For this particular example, we note however that the sub-optimality
gap can be made arbitrarily large by adequately changing
the values for $C(3,1), q(1|3,1)$ and $q(2|3,1)$. These parameters do
not affect the performance of policies that never activate state 3
(including the priority policies in Table~\ref{table_pol}), but do
influence the performance of prio123 and prio213. By making the cost of
being active in state 3, $C(3,1)$, larger, and the transition rates
when being active in state 3 smaller, the performance of these policies
degrades.
As an example, in Figure~\ref{fig:NM_nonindexright} we plot the
sub-optimality gaps when taking $C(3,1)=5$ and
$q(1|3,1)=q(2|3,1)=0.001$ and we observe larger optimality gaps.
Furthermore, note that for a fixed $X(0)$, the gap will grow linearly
in $C(3,1)$, and hence can be made arbitrarily large.

\section{Conclusion and further research}
\label{sec:concl}



In this paper, we studied the general multi-class restless-bandit
problem for both the setting of a fixed population of bandits as well
as a dynamic population of bandits.
Using linear-programming techniques, the paper provided a unified
approach to derive a set of asymptotically optimal priority policies,
$\Pi^*$, which does not rely on indexability of the system.
Under the indexability assumption, Whittle's index policy was shown to
be inside this class.
This is one of the first works that proposes heuristics for \emph
{nonindexable} settings.
As future work, it would therefore be interesting to further understand
their performance outside the asymptotic regime.

The global attractor property is crucial in order to prove asymptotic
optimality of the priority policies $\Pi^*$, as explained in
Section~\ref{sec:global}. Finding sufficient conditions under which the
global attractor property holds for policies in $\Pi^*$ is therefore
important on its own. Another interesting research thread is to
characterize asymptotic optimal policies for models that do not satisfy
the global attractor property, as discussed in Section~\ref{sec:global}.

In addition, it would be interesting to investigate whether
Condition \ref{cond:twee} holds in greater generality for restless
bandit problems.
For example, Condition \ref{cond:twee}(a) concerns stability of the
system under a strict priority policy resulting from the fluid analysis.
In general, care has to be taken when applying a fluid optimal control
directly to the stochastic system, as they might not succeed in making
the system stable; see, for example, \cite{RS92,VN09}.
We believe though that the set $\Pi^*$ contains policies that do
provide a stable system, however, this is a subject for future research.
As an example, we refer to \cite{AEJV13} where a restless bandit
problem was studied that modeled a system with state-dependent capacity.
In that problem, certain priority policies (e.g., the myopic $c\mu$
rule, which is not in $\Pi^*$) yield an unstable system, while other
priority policies, including Whittle's index policy, keep the system stable.

Another interesting research avenue would be to extend this paper to
the general setting of multi-actions. That is, in each state one can
choose from $A_{k}(j)$ different actions, given the constraint $\sum_{k=1}^K\sum_{j=1}^{J_k}\sum_{a=1}^{A_{k}(j)} w^a_{k}(j)
X^a_{j,k}(t)\leq\alpha$, with $w_{k}^a(j)\geq0$ the weight of
action $a$. This paper discussed the case of $w_{k}^a(j)=1$, while
in \cite{HG15} asymptotic optimality of an index policy has been
investigated for $w_{k}^a(j) = a$.

\begin{appendix}\label{app}
\section{Proof of Lemma \texorpdfstring{\protect\ref{15151}}{4.1}}\label{appa}
Set $X_k(0)=x_k(0)$.
Let $\pi$ be a policy for which a unique invariant distribution exists
having finite first moment.
Stability of policy $\pi$ implies rate-stability, that is,
%
\begin{equation}
\label{eq:first} \lim_{t\to\infty} \frac{X_{j,k}^\pi(t)}{t}=0 \qquad\mbox{for
all $j,k$.}
\end{equation}

Note that $\int_{0}^t X_{j,k}^{\pi,a}(s)\,\mathrm{d}s$ is the total
aggregated amount of time spent on action $a$ on class-$k$ bandits in
state $j$ during the interval $(0,t]$.
Hence, we can write the following sample-path construction of the
process $X_{j,k}^\pi(t)$:
%
\begin{eqnarray}
\label{eq:rep1} X_{j,k}^{\pi}(t)
&=&
X_{j,k}^{\pi}(0) + N^{\lambda_k
p_k(j)}(t) + \sum
_{a=0}^{1} \sum_{i=1, i\neq j}^{J_k}
N^{q_{k}(j|i,a)} \biggl(\int_{0}^t
X_{i,k}^{\pi,a}(s)\,\mathrm {d}s \biggr)
\nonumber
\\[-8pt]
\\[-8pt]
\nonumber
&&{} - \sum_{a=0}^{1}\sum
_{i=0, i\neq j}^{J_k} N^{q_k(i| j,a)} \biggl(\int
_{0}^t X_{j,k}^{\pi,a}(s)
\,\mathrm{d}s \biggr),
\end{eqnarray}
where $N^{\lambda_k p_k(j)}(t)$ and $N^{q_{k}(j|i,a)}(t)$ are
independent Poisson processes having as rates $\lambda_k p_k(j)$ and
$q_k(j|i,a)$, respectively, $i,j=1,\ldots, J_k$, $k=1,\ldots, K$, $a=0,1$.
By the ergodic theorem \cite{cinclar75}, we obtain that $\frac{1}{t
}\int_0^{t } X_{j,k}^{\pi, a}(s)\,\mathrm{d}s$ converges to the mean,
denoted by $\overline X^{\pi,a}_{j,k}<\infty$, for all $j,k,a$.
Hence, when dividing both sides in \eqref{eq:rep1} by~$t$, using that
$N^\theta(at)/t\to a\theta$ as $t\to\infty$, and together
with \eqref
{eq:first}, we obtain that
\[
0= \lambda_k p_k(j) + \sum
_{a=0}^{1} \sum_{i=1, i\neq j}^{J_k}
q_{k}(j|i,a) \overline X_{i,k}^{\pi,a} - \sum
_{a=0}^{1} \sum
_{i=0,i\neq j}^{J_k} \overline X_{j,k}^{\pi,a}
q_k(i|j,a) \qquad\mbox{a.s.},
\]
that is, $\overline X^\pi$ satisfies equation \eqref{eq:dif0}. By
definition, $\overline X^\pi$ satisfies $\sum_{k,j } \overline
X_{j,k}^{\pi,1} \leq\alpha$, $\overline X_{j,k}^{\pi,a} \geq0$ and if
$\lambda_k=0$, then $\sum_{j=1}^{J_k}\sum_{a=0}^1 \overline
X_{j,k}^{\pi
, a} = x_k(0)$. Hence, $\overline X^\pi$ is a feasible solution of (LP).

Since the feasible set is nonempty and the objective is to minimize
the cost, the optimal value satisfies $v^*(x(0))<
\infty$.

\section{Proof of Lemma \texorpdfstring{\protect\ref{25252}}{4.3}}\label{appb}
By Fatou's lemma, we have
\[
V^\pi_-(x) \geq\ME_x \Biggl(\liminf
_{T\to\infty} \frac{1}{T}\int_{0}^T
\sum_{k=1}^K \sum
_{j=1}^{J_k}\sum_{a=0}^{1}
C_{k}(j,a) X_{j,k}^{\pi
,a}(t) \,\mathrm{d}t \Biggr).
\]
Hence, it is sufficient to prove that
%
\begin{equation}
\label{eq:bab} \liminf_{T\to\infty} \frac{1}{T}\int
_{0}^T \sum_{k=1}^K
\sum_{j=1}^{J_k}\sum
_{a=0}^{1} C_{k}(j,a)
X_{j,k}^{\pi,a}(t) \,\mathrm{d}t \geq v^*(x)\qquad \mbox{almost
surely},
\end{equation}
with $X(0)=x$.

Consider a fixed realization $\omega$ of the process.
We note that equation \eqref{eq:bab} is trivially true if $\liminf_{T\to\infty} \frac{1}{T}\int_{0}^T \sum_{k=1}^K \sum_{j=1}^{J_k}\sum_{a=0}^{1} C_{k}(j,a) X_{j,k}^{\pi,a}(t)\,\mathrm{d}t=\infty$, since
$v^*(x)<\infty$ (see Lemma~\ref{15151}).
Hence, it remains to be verified that \eqref{eq:bab} holds when the LHS
of \eqref{eq:bab} is finite.

First, assume either a fixed population of bandits, or a dynamic
population of bandits under a stable policy $\pi$.
Since the LHS of \eqref{eq:bab} is finite, we can consider the
subsequence $t_n$ corresponding to the liminf sequence.
For a fixed population of bandits, we have $\frac{1}{T}\int_{0}^T
X_{j,k}^{\pi,a}(t) \,\mathrm{d}t\leq X_{k}(0)=x_k$. Hence, there is a
subsequence $t_{n_l}$ of $t_n$ such that $\frac{1}{t_{n_l}}\int_{0}^{t_{n_l}} X_{j,k}^{\pi,a}(t) \,\mathrm{d}t$ converges to a constant
$\overline X^{\pi,a}_{j,k}$, for all $j,k,a$. In the case of a dynamic
population, given the policy $\pi$ is stable, we have by the ergodicity
theorem \cite{cinclar75} that $\frac{1}{T}\int_0^{T} X_{j,k}^{\pi,
a}(t)\,\mathrm{d}t$ converges to the mean, here denoted by $\overline
X^{\pi,a}_{j,k}$.

In addition, it holds that $\lim_{t\to\infty} X^{\pi,a}_{j,k}(t)/t =0$,
for all $j,k,a$. For the fixed population, this follows since $\lim_{t\to\infty} X^{\pi,a}_{j,k}(t)/t\leq\lim_{t\to\infty} X_{k}(0)/t
=0$, and for the dynamic population this follows since any stable
policy is rate stable.

When studying \eqref{eq:rep1} in the point $t_{n_l}$, dividing both
sides by $t_{n_l}$ and using that $N^\theta(t)/t\to\theta$ as $t\to
\infty$, we can now conclude that
$0 = \lambda_k p_k(j) + \sum_{a=0}^{1} \sum_{i=1}^{J_k} q_{k}(j|i,a)
\overline X_{i,k}^{\pi,a}$.
By \eqref{eq:M} we also have that $\sum_{k=1}^K\sum_{j=1}^{J_k}\overline X_{j,k}^{\pi,1} \leq\alpha$. In addition, if
$\lambda_k=0$, then $\sum_{j,a} \overline X_{j,k}^{\pi, a} =
X_k(0)=x_k$.\vspace*{1pt} Hence $\overline X^\pi$ is a feasible solution of (LP)
with $x(0)=x$.
We conclude that
\begin{eqnarray*}
\label{eq:Valpha} &&\liminf_{T\to\infty} \frac{1}{T}\int
_{0}^T \sum_{k=1}^K
\sum_{j=1}^{J_k}\sum
_{a=0}^{1} C_{k}(j,a)
X_{j,k}^{\pi,a}(t) \,\mathrm{d}t \\
&&\qquad= \lim_{l\to\infty}
\frac{1}{ t_{n_l} }\int_{0}^{t_{n_l}} \sum
_{k=1}^K \sum_{j=1}^{J_k}
\sum_{a=0}^{1} C_{k}(j,a)
X_{j,k}^{\pi,a}(t)\,\mathrm{d}t
\\
&&\qquad= \sum_{k=1}^K \sum
_{j=1}^{J_k}\sum_{a=0}^{1}
C_{k}(j,a) \overline X_{j,k}^{\pi,a} \geq v^*(x),
\end{eqnarray*}
which proves $V_-^\pi(x)\geq v^*(x)$.

We now consider a dynamic population of bandits and take $\pi$ to be
rate-stable. In addition, assume $C_{k}(j,a)>0$, for all $j,k,a$.
Again we consider the subsequence $t_n$ corresponding to the liminf
sequence of \eqref{eq:bab}. So
%
\begin{equation}
\label{eq:sumC} \lim_{n\to\infty} \frac{1}{t_n}\int
_{0}^{t_n} \sum_{k=1}^K
\sum_{j=1}^{J_k}\sum
_{a=0}^{1} C_{k}(j,a)
X_{j,k}^{\pi
,a}(t)\,\mathrm{d}t < \infty.
\end{equation}
Since $C_{k}(j,a)>0$, this implies that the sequence $\frac
{1}{t_n}\int_0^{t_n} X_{j,k}^{\pi, a}(t)\,\mathrm{d}t$ is bounded, for all $j,k,a$.
By the Bolzano--Weierstrass theorem, there exists a subsubsequence
$t_{n_l}$ of $t_n$ and values $\overline X^{\pi,a}_{j,k}$'s such that
$\lim_{l\to\infty}\frac{1}{t_{n_l}} \int_{0}^{t_{n_l}}X^{\pi
,a}_{j,k}(t) \,\mathrm{d}t=\overline X^{\pi,a}_{j,k}$, for all $j,k,a$.
In addition, by rate stability we have that $\lim_{t\to\infty}
X^{\pi
,a}_{j,k}(t)/t =0$, a.s., for all $j,k,a$.
The proof follows now in the same way as above.

The proof in the case of mean-rate stability goes along similar lines
as that for rate stability and is therefore not included here.
%

\section{Proof of Proposition \texorpdfstring{\protect\ref{35353}}{4.13}}\label{appc}
Consider an arbitrary priority policy $\pi$ for which $X^{r,\pi}(t)$ is
irreducible. We first prove stability and then show the tightness and
uniform integrability.

\textit{Stability}:
The Markov process $X^{r,\pi}(t)$ has unbounded transition rates,
however, it does not die in finite time (upward jumps are of the order
1). Hence, once we prove the multi-step drift criterion \cite{MT09,robert03},
we can conclude that there is a unique invariant distribution measure.
The multi-step drift criterion will consist here in proving that there
are $\delta>0$, $T<\infty$, $d>0$ and a stopping time $\tau$, such that
$\ME_x(\tau)\leq T$ for all $x$ and
\[
\ME_x\Biggl(\sum_{k=1}^K
\sum_{j=1}^{J_k} X^{r,\pi,0}_k(
\tau)\Biggr)- \sum_{k=1}^K \sum
_{j=1}^{J_k} x^0_{j,k}\leq-
\delta,
\]
for any $x\in D^c$, with $D:=\{x: \sum_{k=1}^K\sum_{j=1}^{J_k}
x^0_{j,k} \leq d\}$. In other words, for any initial state $x$ outside
the compact set $D$, there is a negative drift (lower bounded by
$-\delta$) toward the set $D$.

We define the stopping time $\tau$ as the first moment that an active
bandit is made passive. Hence, during the interval $[0,\tau]$ the
collection of passive bandits does not change.

First, assume there exists an $x$ such that $\ME_x(\tau)=\infty$. This
implies that when starting in state $x$, the collection of passive and
active bandits remains fixed. Hence, each passive class-$k$ bandit
evolves according to the transition rates $q_k(j|i,0)$. The number of
passive class-$k$ bandits is therefore equivalent to that in an
$M/G/\infty$ queue with arrival rate $\lambda_k r$ and phase-type
distributed service requirements as described by the transitions of a
passive class-$k$ bandit. We note that the $M/G/\infty$ queue is stable.
By irreducibility, for any starting point, the process will be in
state $x$ after a finite expected amount of time, hence, stability follows.

We now assume $\ME_x(\tau)<\infty$, for all $x$. Since there is a
finite number of states $J_k<\infty$ and the state transitions are
exponential, it follows directly that there exists a $T<\infty$ such
that $\ME_x(\tau)<T$, for all $x$.
Note that the passive bandits behave independently during the interval
$[0,\tau]$.
The probability that a passive bandit departs in the interval $[0,\tau
]$ can be lower bounded by $p_0$ with $p_0>0$.
This follows from the assumption that state $0$ is positive recurrent
under the policy that always keeps the class-$k$ bandit passive.
Hence, the mean number of passive bandits that leave during the
interval $[0,\tau]$ is larger than or equal to $p_0 \sum_{k=1}^{K}
\sum_{j=1}^{J_k} x^0_{j,k}$.
We therefore have as mean drift
\begin{eqnarray*}
&&\ME_x\Biggl(\sum_{k=1}^K
\sum_{j=1}^{J_k} X^{r,\pi,0}_k(
\tau)\Biggr)- \sum_{k=1}^K \sum
_{j=1}^{J_k} x^0_{j,k}
\\
&&\qquad\leq \lambda r \ME_x(\tau)+1-p_0 \sum
_{k}^{K} \sum_{j=1}^{J_k}
x^0_{j,k}<\lambda r T+1-p_0 d,
\end{eqnarray*}
for all $x\in D^c$.
The $+1$ in the mean drift is due to the active bandit that becomes
passive at time $\tau$. Choosing $d=(\lambda r T+1+\delta)/p_0$, we
conclude that $\ME_x(\sum_{k=1}^K \sum_{j=1}^{J_k} X^{r,\pi
,0}_k(\tau
))- \sum_{k=1}^K \sum_{j=1}^{J_k} x^0_{j,k}\leq- \delta$.
Hence, by the multi-step drift criterion we obtain that there is a
unique invariant probability distribution for the process $X^{r,\pi
}(t)$, for any $r$. Recall that we denote this distribution by
$p^{r,\pi}$.

\textit{Tightness and uniform integrability}:
In order to prove tightness and uniform integrability, we will define a
process that serves as a stochastic upper bound on $X_{k}^{r,\pi,0}(t)$.
First, note that $\max_{i} q_k(j|i,1)\alpha$ is the maximum rate at
which active bandits go to state $j$. Hence, $\overline\lambda_{k} :=
\lambda_k + \sum_{j=1}^{J_k}\max_{i} q_k(j|i,1)\alpha$ is an upper
bound on the arrival rate of new passive class-$k$ bandits. For the
upper bound process, we assume that once a bandit is passive, it will
never be made active again. Hence, the time such a passive bandit stays
in the system can be described by the state transitions of a passive
class-$k$ bandit. We define $B_k$ as the distribution described by the
state transition rates $q_k(j|i,0)$, with a certain initial probability
$\tilde p_0$. Choosing an appropriate value for $\tilde p_0$, the $B_k$
describes the time a passive class-$k$ bandit stays in the system. Let
$Y^r_k(t)$ be the number of customers in a $M/G/\infty$ queue with
arrival rate $\overline\lambda_k$ and service requirement $B_k$. This
process is an upper bound on $X_{k}^{r,\pi,0}(t)$.


The stationary distribution of the process $\{Y_{k}^r(t)\}$ is given by
a Poisson distribution with parameter $\overline\lambda_{k} r \ME
(B_{k})$ \cite{robert03}.
It can be checked that this distribution converges to the Dirac measure
in the point $\overline\lambda_{k} \ME(B_{ k})$, as $r\to\infty$.
By Prohorov's theorem, it then follows that the family $\{Y_{k}^r/r\}$
is tight \cite{robert03}.
Furthermore, since $\ME(Y_{k}^r/r) = \overline\lambda_{k} \ME(B_{k})$
and $ \ME( \lim_{r\to\infty} Y_{k}^r/r) =\overline\lambda_{k} \ME
(B_{k})$, a.s., we obtain from \cite{B99}, Theorem~3.6, that the family
$\{Y_{k}^r/r\}$ is uniform integrable.

At most $\alpha$ bandits are active, hence $\sum_{ k} Y_{
k}^r(t)/r+\alpha$ represents a stochastic upper bound on the queue
length process $\sum_{k=1}^K X_k^{r,\pi}(t)/r$. This implies that the
family $\{p^{r,\pi}\}$ is tight and uniform integrable as well.

\section{Proof of Proposition \texorpdfstring{\protect\ref{45454}}{4.14}}\label{appd}
We denote by $S^{\pi^*}_{k}(j)$ the set of all combinations $(i,l),
i=1,\ldots, J_l, l=1,\ldots, K$, such that class-$l$ bandits in
state $i$ have higher priority than class-$k$ bandits in state $j$
under policy $\pi^*$, and $I^{\pi^*}$ is the set of all states that
will never be made active under policy $\pi^*$.
The transition rates of the process $X^{r, \pi^*}(t)/r$ are then
defined as follows:
%
\begin{eqnarray}
x &\to& x + \frac{e_{j,k}}{r}\qquad \mbox{at rate } r \lambda_k
p_k({j}), k=1,\ldots, K, j=1,\ldots, J_k,
\label{eq:eq1}
\\[-2pt]
\label
{eq:eq2}
x &\to& x - \frac{e_{j,k}}{r}
\nonumber
\\[-14pt]
\\[-12pt]
\eqntext{\mbox{at rate } r\displaystyle  \sum
_{a=0}^1 x_{j,k}^{a}
q_k(0|j,a), k=1,\ldots, K, j=1,\ldots, J_k,}
\\[-2pt]
\label{eq:eq3}
x &\to& x - \frac{e_{j,k}}{r}+ \frac{e_{i,k}}{r}
\nonumber
\\[-12pt]
\\[-12pt]
\eqntext{\mbox{at rate } r \displaystyle \sum
_{a=0}^1 x_{j,k}^{a}
q_k(i|j,a), k=1,\ldots, K, i,j=1,\ldots, J_k, i\neq j,}\vspace*{-2pt}
\end{eqnarray}
where $x^{1}_{j,k}=\min ((\alpha- \sum_{(i,l) \in S^{\pi^*}_{
k}(j)} x_{i,l})^+, x_{j,k}  )$, if $(j,k)\notin I^{\pi^*}$, and
$x^{1}_{j,k}=0$ otherwise, $x^{0}_{j,k}= x_{j,k}- x^{1}_{j,k}$, and
$e_{j,k}$ is a vector composed of all zeros except for component
$(j,k)$ which is one.

From \eqref{eq:eq1}--\eqref{eq:eq3}, it follows that there exists a
continuous function $b_l(x)$, with $l\in{\cal{L}}$ and ${\cal{L}}$
composed of a finite number of vectors in $\mathbb{N}^{\sum_k J_k}$,
such that the transition rates of the process $x^{r,\pi^*}(t)$ from $x$
to $x+ l/r$ have the form $r b_l(x)$.
Hence, the process $X^{r, \pi^*}_{j,k}(t)/r$ belongs to the family of
density dependent population processes as defined in \cite{EK86}, Chapter~11.

Note that the process $x^{\pi^*}(t)$ as defined in \eqref{eq:ode} can
equivalently be written as $\frac{\mathrm{d}x^{\pi^*}(t)}{\mathrm
{d}t}=F(x^{\pi^*}(t))$, with $F(x^*)=\sum_{l\in{\cal{L}}} lb_l(x^*)$,
where $F(\cdot)$ is Lipschitz continuous.
From Condition \ref{cond:better}, we have that $x^*$ is the unique
global attractor of $x^{\pi^*}(t)$.

Together with the fact that the family $\{p^{r,\pi^*}\}$ is tight, we
then obtain from \cite{GG10}, Theorem~4, that $p^{r,\pi^*}(x)$ converges
to the Dirac measure in $x^*$, the global attractor of $x^{\pi^*}(t)$.
Hence, we can write
\begin{eqnarray*}
\lim_{r\to\infty} V_+^{r,\pi^*}(x) &=& \sum
_{k=1}^K\sum_{j=1}^{J_k}
\sum_{a=0}^1 \lim_{r\to\infty}
\sum_{x} p^{r,\pi^*}(x) C_{k}(j,a)
x^{a}_{j,k} \\[-3pt]
&=& \sum_{k=1}^K
\sum_{j=1}^{J_k} \sum
_{a=0}^1 C_{k}(j,a)
x^{*,a}_{j,k} = v^*(x),
\end{eqnarray*}
where the first step follows from the ergodicity theorem \cite
{Tijms2003,cinclar75} (applicable since the first moment of $p^{r,\pi
^*}$ is finite), the second step (interchange of limit and summation)
follows from uniform integrability of $\{p^{r,\pi^*}\}$ and the fact
that $p^{r,\pi^*}$ converges to the Dirac measure in $x^*$, and the
last step follows since $x^*$ is an optimal solution of (LP).


We conclude the proof by noting that $v^*(x)$ is a lower bound on the
steady-state cost, as shown in Lemma~\ref{25252}.

\section{Proof of Proposition \texorpdfstring{\protect\ref{56565}}{5.6}}\label{appe}
Recall that the relaxed optimization problem for $f=\mathrm{av}$ consists in
finding a stationary and Markovian policy that minimizes
%
\begin{equation}
\label{eq:ApB1} \lim_{T\to\infty} \frac{1}{T}
\mathbb{E}_x \Biggl(\int_0^T \sum
_{k=1}^K \sum
_{j=1}^{J_k} \sum_{a=0}^{1}
C_{k}(j,a) X_{j,k}^{\pi,a}(t)\,\mathrm {d}t \Biggr),
\end{equation}
under the relaxed constraint
%
\begin{equation}
\label{eq:ApB1b} \lim_{T\to\infty} \frac{1}{T}
\mathbb{E}_x \Biggl(\int_0^T \sum
_{k=1}^K \sum
_{j=1}^{J_k} X_{j,k}^{\pi,1}(t)
\,\mathrm{d}t \Biggr)\leq\alpha.
\end{equation}

For a given policy $\pi$, we denote by $x_{j,k}^{\pi,a}$ the
(stationary) state-action frequencies, that is, the average fraction of
time the class-$k$ bandit is in state $j$ and action $a$ is chosen.
Assumption~\ref{as:unichain} implies that these frequencies exist and
satisfy the balance equations, that is, they satisfy
\[
0 = \sum_{a=0}^1 \sum
_{i=0, i\neq j}^{J_k}q_k(i|j,a)
x^{\pi,a}_{j,k} - \sum_{a=0}^1
\sum_{i=1, i\neq j}^{J_k} q_k(j|i,a)
x^{\pi
,a}_{i,k}\qquad\forall j, 
\]
or, by definition of $q_k(j|j,a)=-\sum_{i=0, i\neq j}^{J_k}
q_k(i|j,a)$, this can be written as
\[
0 = \sum_{a=0}^1 \sum
_{i=1 }^{J_k} q_k(j|i,a)
x^{\pi,a}_{i,k} \qquad\forall j. 
\]

We will restrict ourselves to the class of policies that are symmetric
for bandits in the same class. We can do this without loss of
generality, since an optimal solution of the relaxed problem, given by
Whittle's indices, is symmetric.
Having $X_k(0)$ bandits in class $k$, equations \eqref{eq:ApB1}
and \eqref{eq:ApB1b} can now equivalently be written as
\[
\sum_{k=1}^K X_k(0) \sum
_{j=1}^{J_k} \bigl( C_{k}(j,0)
x^{\pi
,0}_{j,k} + C_{k}(j,1) x^{\pi,1}_{j,k}
\bigr)\quad \mbox{and}\quad \sum_{k=1}^KX_k(0)
\sum_{j=1}^{J_k}x^{\pi,1}_{j,k}
\leq\alpha,
\]
respectively.

The relaxed optimization problem can now be formulated as the following
linear program (D):
\begin{eqnarray} \label{eq:ApB2}
\mathrm{(D)} \quad &&\min_x \sum_{k=1}^K
X_k(0) \sum_{j=1}^{J_k} \bigl(
C_{k}(j,0) x^{0}_{j,k} + C_k(j,1)
x^1_{j,k} \bigr)
\nonumber
\\
&&\mbox{s.t.}\qquad 0 = \sum_{a=0}^1 \sum
_{i=1 }^{J_k} q_k(j|i,a)
x^a_{i,k}\qquad \forall j,k,
\nonumber
\\[-8pt]
\\[-8pt]
\nonumber
&& \hphantom{\mbox{s.t.}\qquad}\sum_{k=1}^K X_k(0) \sum
_{j=1}^{J_k}x^1_{j,k}
\leq\alpha,
\\
&& \hphantom{\mbox{s.t.}\qquad}\sum_{j=1}^{J_k} \sum
_{a=0}^1 x^a_{j,k} =1\qquad
\forall k,\qquad x^a_{j,k}\geq0\qquad \forall k,j,a.
\nonumber
\end{eqnarray}
We have that for any feasible solution $(x^{a}_{j,k})$ of (D) there is
a stationary policy $\pi$ such that the state-action frequencies
$x^{\pi
,a}_{j,k}$ coincide with the value of the feasible solution
$x^{a}_{j,k}$ \cite{P94}, Theorem~8.8.2(b).
Hence, for any optimal (symmetric) policy $\pi^*$ of the relaxed
optimization problem, the state-action frequencies $x_{j,k}^{\pi^*,a}$
provide an optimal solution of (D). We further note that $(x^{\pi
^*}_{j,k}X_k(0))$ is an optimal solution of (LP) with $x(0)=X(0)$.

We assume the restless bandit problem is indexable. Hence, an optimal
policy of the relaxed optimization problem is described in Section~\ref{sec:rel}, and will be denoted here by $\tilde\pi^*$.
We recall that policy $\tilde\pi^*$ is described by a value $\nu
^*\geq
0$ and is such that a class-$k$ bandit in state $j$ is made active if
$\nu_{k}^{\mathrm{av}}(j)>\nu^*$ and is kept passive if $\nu_{k}^{\mathrm{av}}(j)<\nu^*$.
Hence, the state-action frequencies under $\tilde\pi^*$ satisfy
%
\begin{eqnarray}\label{eq:jee}
x^{\tilde\pi^*,0}_{j,k}&=&0\qquad \mbox{when } \nu_{k}^{\mathrm{av}}(j)>
\nu^*,
\nonumber
\\[-8pt]
\\[-8pt]
\nonumber
x^{\tilde\pi^*,1}_{j,k}&=&0\qquad \mbox{when } \nu_{k}^{\mathrm{av}}(j)
< \nu^*.
\end{eqnarray}
By definition of policy $\tilde\pi^*$, for states $(\hat j,\hat k)$
with $\nu^{\mathrm{av}}_{ \hat k}(\hat j,)=\nu^*$ a class-$\hat k$ bandit in
state $\hat j$ is made active with a certain probability, hence
$x^{\tilde\pi^*,0}_{\hat j,\hat k}\geq0$ and $x^{\tilde\pi
^*,1}_{\hat j,\hat k}\geq0$.

Since Whittle's index policy gives priority to bandits having highest
index value, we directly obtain that Whittle's index policy $\nu^{\mathrm{av}}$
satisfies points 1 and 2 of Definition~\ref{def:Pix} when setting
$x^*=(x_{j,k}^{\tilde\pi^*}X_k(0))$.\vspace*{-1pt}
We now treat point 3 of Definition~\ref{def:Pix}: Assume $\sum_{k=1}^K\sum_{j=1}^{J_k} x^{\tilde\pi^*,1}_{j,k}X_k(0) < \alpha$.
Hence, under the optimal policy, on average, strictly less than $\alpha
$ bandits are made active. This implies that the remaining fraction of
the time the policy makes dummy bandits in state $B$ active. Hence,
$\nu
_B^{\mathrm{av}}\geq\nu^*$. Since $\nu^*\geq0$ and $\nu_B^{\mathrm{av}}=0$, we
necessarily have $\nu^*=0$.
A policy satisfies point 3 of Definition~\ref{def:Pix} if it never
makes a class-$k$ bandit in state $j$ active that satisfies
%
\begin{equation}
\label{eq:moe} x_{j,k}^{\tilde\pi^*,1}=0\quad \mbox{and}\quad
x_{j,k}^{\tilde
\pi^*,0}>0.
\end{equation}
From \eqref{eq:jee} (with $\nu^*=0$), we obtain that \eqref{eq:moe}
implies $\nu_{k}^{\mathrm{av}}(j)\leq0$.
By definition of Whittle's index policy, a bandit in a state such that
$\nu_{k}^{\mathrm{av}}(j)\leq0$ will never be made active, hence point 3 is
satisfied. We therefore conclude that Whittle's index policy $\nu^{\mathrm{av}}$
is included in the set of priority policies $\Pi(x^* )\subset\Pi^*$,
with $x^*=(x_{j,k}^{\tilde\pi^*}X_k(0))$. 

\section{Proof of Proposition \texorpdfstring{\protect\ref{75757}}{5.9}}\label{appf}
Let $\beta\leq\overline\beta$ and $\beta>0$.
Whittle's index $\nu_{k}^\beta(j)$ results from solving the following
problem for a class-$k$ bandit:
%
\begin{equation}
\label{eq:timeav2} \min_{\pi} \mathbb{E}_x \biggl(
\int_{0}^\infty\mathrm{e}^{-\beta
t}
\bigl(C_k\bigl(J_k(t), A_k^\pi(t)
\bigr) + \nu\mathbf{1}_{(A_k^\pi(t)=1)} \bigr)\, \mathrm {d}t \biggr),
\end{equation}
see \eqref{eq:subproblem}, where $A_k^\pi(t)\in\{0,1\}$ and $J_k(t)$
denotes the state of the class-$k$ bandit.
This is a continuous-time discounted Markov decision problem in a
finite state space. After uniformization (\cite{GHP06}, Remark~3.1,
\cite{P94}, Section~11.5.2), this is equivalent to a \emph
{discrete-time} discounted Markov decision problem with discount factor
$\tilde\beta=\frac{\overline q}{\beta+\overline q}$, cost function
$\tilde C_k(j,a) =\frac{C_{k}(j,a)+ \nu\mathbf{1}_{(a=1)}}{\beta+
\overline q}$, and transition probabilities $\tilde p^a_k(i,j) = \frac
{q_k(j|i,a)}{\overline q} + \mathbf{1}_{(i=j)}$ [recall that
$q_k(i|i,a) =\break -\sum_{j=0, i\neq j}^{J_k} q^a_k(i,j)$], where
$\overline
q := \max_{i,k,a} -q_k(i|i,a)<\infty$.
In LP formulation the discrete-time MDP for the class-$k$ bandit is
then as follows (see \cite{P94}, Section~6.9):
\begin{eqnarray}
&&\max_{v} \sum_{j=1}^{J_k}
\gamma_{j,k} v(j)
\nonumber
\\
&&\mbox{s.t.}\qquad v(i)-\tilde\beta\sum_{j=0}^{J_k}
\tilde p_k^a(i,j) v(j) \leq\tilde C_{k}(i,a)\qquad
\forall i=1,\ldots, J_k, a=0,1,
\nonumber
\end{eqnarray}
with $\gamma_{j,k}>0$ arbitrary. In fact, we will make the choice
$\gamma_{j,k}=\lambda_k (p_{0j}^k + \varepsilon)$, with $\varepsilon>0$.
The dual of the above LP is
%
\begin{eqnarray}\label{eq:bc} \label{eq:bcd}
\mathrm{\bigl(D_k(\beta,\varepsilon)\bigr)} \quad &&\min_x
\sum_{j=1}^{J_k}\frac{ C_k(j,0)
x^0_{j,k} + C_{k}(j,1) x^1_{j,k} +\nu x^1_{j,k}}{\beta+\overline
q}
\nonumber
\\
&&\mbox{s.t.}\qquad 0= \lambda_k \bigl(p_k({0,j})+
\varepsilon\bigr)
\nonumber
\\[-8pt]
\\[-8pt]
\nonumber
&&\hspace*{52pt}{}+ \sum_{a=0}^1 \sum
_{i=1 }^{J_k} \frac{q_k({j|i,a})}{\beta+ \overline q}
x^a_{i,k} - \frac{\beta}{\beta+\overline{q}} \sum
_{a=0}^1 x^a_{j,k}\quad \forall
j,
\\
&&\hphantom{\mbox{s.t.}\qquad}  x^a_{j,k}\geq0 \qquad\forall j, a.
\nonumber
\end{eqnarray}

As stated in Section~\ref{subsec:relax}, indexability implies that an
optimal policy for the subproblem \eqref{eq:timeav2} is described by a
priority ordering according to the indices $\nu_{k}^\beta(j)$: an
optimal action in state $j$ is $a =1$ if $\nu_{k}^\beta(j)> \nu$ and $a
=0$ if $\nu_{k}^\beta(j)< \nu$.
Recall that a class-$k$ bandit is indexable for each $\beta_l$
(subsequence can depend on the class). Hence, by \cite{P94}, Theorem~6.9.4, this implies that there exists an optimal solution to
[$D_k(\beta_l,\varepsilon)$], denoted by $x_k^*(\beta_l, \varepsilon)$,
such that
\begin{eqnarray*}
 x^{*,0}_{j,k}(\beta_l, \varepsilon)&=&0 \qquad\mbox{when } \nu _{k}^{\beta
_l}(j)> \nu,
\\
 x^{*,1}_{j,k}(\beta_l, \varepsilon)&=&0\qquad \mbox{when } \nu _{k}^{\beta
_l}(j) < \nu.
\end{eqnarray*}
Since $\lim_{l\to\infty} \nu_{k}^{\beta_l}(j)=\nu_{k}^{\mathrm{lim}}(j)$, we
obtain that there exists an $L(\nu)$ such that for all $l>L(\nu)$ it
holds that
%
\begin{eqnarray}
x^{*,0}_{j,k}(\beta_l, \varepsilon)&=&0\qquad \mbox{when } \nu_{k}^{\mathrm{lim}}(j)> \nu, \label{eq:prop1_d}
\\
x^{*,1}_{j,k}(\beta_l, \varepsilon)&=&0 \qquad\mbox{when } \nu_{k}^{\mathrm{lim}}(j) <\nu.\label{eq:prop1_de}
\end{eqnarray}


By change of variable $\tilde x_{j,k}^{a}= x_{j,k}^{a}/(\beta+
\overline q)$ we obtain that $\tilde x_k^{*}(\beta_l, \varepsilon)$
satisfies \eqref{eq:prop1_d} and \eqref{eq:prop1_de} and is an optimal
solution of [$\tilde D_k(\beta_l, \varepsilon)$] defined as
%
\begin{eqnarray} \label{eq:jj}
\bigl(\tilde D_k(\beta,\varepsilon)\bigr) \quad &&\min
_{\tilde x} \sum_{j=1}^{J_k}
\bigl(C_{k}(j,0) \tilde x^0_{j,k} +
C_{k}(j,1) \tilde x^1_{j,k} +\nu\tilde
x^1_{j,k}\bigr)
\nonumber
\\
&&\mbox{s.t.}\qquad 0 = \lambda_k \bigl(p_k(j)+\varepsilon
\bigr)
\nonumber
\\[-8pt]
\\[-8pt]
\nonumber
&&\hspace*{52pt}{} + \sum_{a=0}^1 \sum
_{i=1, i\neq j}^{J_k} q_k({j|i,a}) \tilde
x^a_{i,k} - \beta\sum_{a=0}^1
\tilde x^a_{j,k}\qquad \forall j,
\\
&& \hphantom{\mbox{s.t.}\qquad}\tilde x^a_{j,k}\geq0\qquad\forall j, a.
\nonumber
\end{eqnarray}

By Assumption~\ref{as:bounded2}, we have that the set of optimal
solutions of $(\tilde D_k(0,0))$ is bounded and nonempty when $\nu>0$.
Hence, from \cite{CLP05}, Corollary~1, we obtain that the correspondence
that gives for each $(\beta,\varepsilon)$ the set of optimal solutions of
$(\tilde D_k(\beta,\varepsilon$)) is upper semi-continuous in the point
$(\beta, \varepsilon)=(0,0)$. It is a compact-valued correspondence [after
summing \eqref{eq:jj} over all $j$, we have that $\tilde x_{k}=\lambda
_k(1+\varepsilon J_k)/ \beta$, $\beta>0$]. Hence, it follows that there
exists a sequence $(\beta_{l_n}, \varepsilon_n)$ (with $\beta_{l_n}$ a
subsequence of $\beta_l$ and $\varepsilon_n\to0$) such that $\tilde
x^{*,a}_{j,k}(\beta_{l_n}, \varepsilon_n)\to\tilde x^{*,a}_{j,k}$, as
$n\to\infty$, and with $\tilde x^*_k$ an optimal solution of $(\tilde
D_k(0,0))$.
For a fixed $\nu$, the components of $\tilde x^*_k(\beta_l, \varepsilon)$
that are zero are independent of the exact values for $\varepsilon>0$, and
$l>L(\nu)$; see \eqref{eq:prop1_d} and \eqref{eq:prop1_de}. Hence, the
limit $\tilde x^*_k$, which is an optimal solution of$ (\tilde
D_k(0,0))$, has the same components equal to zero, that is, \eqref
{eq:prop1_d} and \eqref{eq:prop1_de} are satisfied for $\tilde x^*_k$.

Below we will show that there exists a value $\nu^*$ such that there is
a vector $\tilde y^{*}$ that satisfies the following: (i) $\tilde
y_{k}^{*}$ is an optimal solution of $(\tilde D_k(0,0))$, for all $k$,
with $\nu=\nu^*$, (ii) $\tilde y^*$ is an optimal solution of (LP), and
(iii) the Whittle index policy $\nu^{\mathrm{lim}}$ is included in the set $\Pi
(\tilde y^*)\in\Pi^*$. The latter then completes the proof.

In the remainder of the proof, we denote by $\tilde x_k^*(\nu)$ the
above described optimal solution $\tilde x_k^*$ of $(\tilde D_k(0,0))$
for a given value $\nu$.
We have the following properties:
\begin{itemize}
\item\textit{Property}~1:
%
\begin{equation}
\sum_{k=1}^K\sum
_{j=1}^{J_k} \tilde x^{*,1}_{j,k}(
\infty) \leq \alpha. \label{eq:M1}
\end{equation}
This can be seen as follows.
As $\nu\to\infty$, the objective of $(\tilde D_k(0,0))$ is to minimize
$\sum_{j=1}^{J_k} \tilde x^1_{j,k}$.
For any feasible solution $x$ of (LP), $x_k$ is in the feasible set of
$\tilde D_k(0,0)$.
Hence, $\sum_{j=1}^{J_k} \tilde x^{*,1}_{j,k}(\infty) \leq\sum_{j=1}^{J_k} x^{1}_{j,k}$ with $x$ a feasible solution of (LP).
In addition, we have that $\sum_{k=1}^K \sum_{j=1}^{J_k} x^{1 }_{j,k}
\leq\alpha$ with $x$ a feasible solution of (LP). This proves \eqref{eq:M1}.
\item
\textit{Property}~2:
%
\begin{equation}
\label{eq:M2} \sum_{j=1}^{J_k} \tilde
x^{*,1}_{j,k}(\nu) \geq\sum_{j=1}^{J_k}
\tilde x^{*,1}_{j,k}(\tilde\nu)\qquad \mbox{for } \nu<\tilde\nu.
\end{equation}
This can be seen as follows:
By definition, we have
$\sum_{j=1}^{J_k}\sum_{a=0}^1 C_{k}(j,a) \tilde x^{*,a}_{j,k}(\nu) +
\nu\sum_{j=1}^{J_k} \tilde x^{*,1}_{j,k}(\nu) \leq\sum_{j=1}^{J_k}\sum_{a=0}^1 C_{k}(j,a) \tilde x^{*,a}_{j,k}(\tilde\nu) +
\nu\sum_{j=1}^{J_k} \tilde x^{*,1}_{j,k}(\tilde\nu)$ and\break 
$\sum_{j=1}^{J_k}\sum_{a=0}^1 C_{k}(j,a) \tilde x^{*,a}_{j,k}(\tilde
\nu
) +\tilde\nu\sum_{j=1}^{J_k} \tilde x^{*,1}_{j,k}(\tilde\nu)
\leq\sum_{j=1}^{J_k}\sum_{a=0}^1 C_{k}(j,a) \tilde x^{*,a}_{j,k}(\nu)
+ \tilde\nu\sum_{j=1}^{J_k} \tilde x^{*,1}_{j,k}(\nu)$.
Subtracting the latter inequality from the first, we obtain
equation \eqref{eq:M2}.

%
\item
\textit{Property}~3:
%
\begin{equation}
\label{eq:M3} \sum_{j=1}^{J_k} \tilde
x^{*,1}_{j,k}(\nu)<\infty \qquad\mbox{for } \nu>0.
\end{equation}
This follows since by Assumption~\ref{as:bounded2} the set of optimal
solutions of $(\tilde D_k(0,0))$ is bounded for $\nu>0$.
\end{itemize}
%
We define $\overline\alpha:= \sum_{k=1}^K\sum_{j=1}^{J_k} \tilde
x^{*,1}_{j,k}(0)$. Equations \eqref{eq:M1}--\eqref{eq:M3} imply that
there exists a $\nu^*\geq0$ such that
%
\begin{equation}\qquad
\label{eq:ff} \sum_{k=1}^K\sum
_{j=1}^{J_k} \tilde x^{*,1}_{j,k}
\bigl(\bigl(\nu^*\bigr)^-\bigr) \geq \min (\alpha, \overline\alpha) \quad\mbox{and}\quad
\sum_{k=1}^K\sum
_{j=1}^{J_k} \tilde x^{*,1}_{j,k}
\bigl(\bigl(\nu^*\bigr)^+\bigr) \leq\min(\alpha, \overline\alpha).
\end{equation}
From standard LP theory, we know that there exists a $\overline\nu
<\infty$ such that $\tilde x^*_k(\overline\nu)$ is an optimal solution
of $(D_k(0,0))$ for all $\nu\geq\overline\nu$, that is $\tilde
x^*_k(\nu)=\tilde x^*_k(\overline\nu)$ for $\nu\geq\overline\nu$.
Hence, we can take $\nu^*<\infty$.

From \eqref{eq:M2} and \eqref{eq:ff}, we obtain that there exists a
$\tilde y^*=(\tilde y^{*,a}_{j,k})$ with $y^*_{\tilde k}$ being a
convex combination of $\tilde x_{\tilde k}^*((\nu^*)^-)$ and $\tilde
x_{\tilde k}^*((\nu^*)^+)$ and for $k\neq\tilde k$, $\tilde y^*_k$
being equal to either $\tilde x_k^*((\nu^*)^-)$ or $\tilde x^*_k((\nu
^*)^+)$, such that $\sum_{k=1}^K\sum_{j=1}^{J_k} \tilde
y^{*,1}_{j,k}=\min(\alpha, \overline\alpha)$. Note that $\tilde y^*_k$
is still a solution of $(\tilde D_k(0,0))$, for all $k$.
Now, if $\alpha= \min(\overline\alpha, \alpha)$, it follows directly
that $\tilde y^*$ is also an optimal solution of (LP).
If instead $\overline\alpha= \min(\overline\alpha, \alpha)$, then
$\nu^*=0$, and hence $\tilde y^*_k$ is an optimal solution of
$(\tilde
D_k(0,0))$ with $\nu=0$. After summing over $k$, the latter has the
same objective function as (LP). Together with $\sum_{k=1}^K \sum_{j=1}^{J_k} \tilde y_{j,k}^{*,1} = \overline\alpha\leq\alpha$, it
follows that $\tilde y^*$ is also an optimal solution of (LP).

It remains to be proved that the Whittle index policy is included in
the set $\Pi(\tilde y^*)\subset\Pi^*$.
Assume for class $\tilde k$ the states are ordered such that $\nu_{
{\tilde k} }^{\mathrm{lim}}(j_1) \leq\nu_{ {\tilde k} }^{\mathrm{lim}}(j_2)<\cdots\leq
\cdots\leq\nu_{ {\tilde k} }^{\mathrm{lim}}(j_{J_{\tilde k}})$.
From $\nu^*<\infty$ and properties \eqref{eq:prop1_d}--\eqref
{eq:prop1_de} [which hold for $\tilde x^*(\nu)$], we have that there
are $n^*$ and $\tilde n^*$, $n^*\leq\tilde n^*$, such that $\nu
_{{\tilde k}} (j_{n^*})= \cdots = \nu_{{\tilde k} }(j_{\tilde
n^*})=\nu^*$ and
\begin{eqnarray*}
\tilde x_{j_m, {\tilde k} }^{*,1}\bigl(\bigl(\nu^*\bigr)^-\bigr) &=& 0\qquad
\mbox{for all } m=1,\ldots, n^*,
\\
\tilde x_{j_m, {\tilde k} }^{*,0}\bigl(\bigl(\nu^*\bigr)^-\bigr) &=& 0\qquad
\mbox{for all } m=n^*+1,\ldots, J,
\end{eqnarray*}
and
\begin{eqnarray*}
\tilde x_{j_m, {\tilde k} }^{*,1}\bigl(\bigl(\nu^*\bigr)^+\bigr) &=&0\qquad
\mbox{for all } m=1,\ldots,\tilde n^*,
\\
\tilde x_{j_m, {\tilde k} }^{*,0}\bigl(\bigl(\nu^*\bigr)^+\bigr) &=& 0\qquad
\mbox{for all } m=\tilde n^*+1,\ldots, J.
\end{eqnarray*}
%
The vector $\tilde y_{\tilde k}^*$ is a convex combination of $\tilde
x_{\tilde k}^*((\nu^*)^-)$ and $\tilde x_{\tilde k}^*((\nu^*)^+)$,
hence $\tilde y_{j_m, \tilde k}^{*,1}=0$ for all $m\leq n^*$ and
$\tilde y_{j_m, \tilde k}^{*,0}=0$ for all $m\geq\tilde n^*+1$. Hence,
Whittle's index policy $\nu^{\mathrm{lim}}$ satisfies items 1 and 2 of
Definition~\ref{def:Pix} with $x^*=\tilde y^*$.

If $\sum_{k=1}^K \sum_{j=1}^{J_k} \tilde y_{j,k}^{*,1}<\alpha$, then
since $\sum_{k=1}^K \sum_{j=1}^{J_k} \tilde y_{j,k}^{*,1}=\min
(\alpha,
\overline\alpha)$ we have $\bar\alpha< \alpha$, so $\nu^*=0$. This
implies that for any state $(j,k)$ with $\tilde y^{*,1}_{j,k}=0$ and
$\tilde y^{*,0}_{j,k}>0$ it follows from property \eqref{eq:prop1_d}
that $\nu^{\mathrm{lim}}_{k}(j) < (\nu^*)^+=0^+$. Hence, by definition of
Whittle's index policy $\nu^{\mathrm{lim}}$, a bandit in this state will never
be made active, which implies that item 3 in Definition~\ref{def:Pix}
is satisfied for $x^*=\tilde y^*$. It hence follows that Whittle's
index policy $\nu^{\mathrm{lim}}$ is included in the set of priority
policies $\Pi(\tilde y^*)\subset\Pi^*$.
%
%
%

\section{Proof of Lemma \texorpdfstring{\protect\ref{888885}}{8.1}}\label{appg}
For the fixed population, the total number of constraints in (LP) is
$\sum_{k=1}^K J_k +1 +K$.
However, since $\sum_{k=1}^K \lambda_k=0$, one of the constraints
in \eqref{eq:dif0} is redundant for each $k$. Hence, the number of
independent constraints in (LP) is $\sum_{k=1}^K J_k+1$.

Since the feasible set of (LP) is bounded, from standard LP theory
(see \cite{P94}, Theorem D.1a), we obtain that there exists an optimal
basic feasible solution $x^*$ to (LP). Hence, $x^*$ has $\sum_{k=1}^K
J_k+1$ basic terms and all other terms are equal to zero.
If $x_{j,k}^*>0$ for all $j,k$, then for any $j,k$ there is an action
$a$ such that $x_{j,k}^{*,a}=0$, and in at most one combination $(j,k)$
the components $x_{j,k}^{*,a}$ can be positive in both actions. Hence,
$x^*$ satisfies the property in Definition~\ref{def:Pix}.

Otherwise, let $S$ denote the set of pairs $(i, l)$ such that
$x_{i,l}^*=0$. By \eqref{eq:dif0}, if $(j,k)\in S$, then $\sum_{a=0}^1
\sum_{i\neq j} x_{i,k}^{*,a}q_k(j|i,a)=0$. That is,
$x_{i,k}^{*,a}q_k(j|i,a)=0$ for all $i=1,\ldots, J_k$, $a=0,1$, if
$(j,k)\in S$.
Hence, for $(j,k)\notin S$, equation \eqref{eq:dif0} in the point $x^*$
can be rewritten as
\[
0 = \sum_{a=0}^{1} \sum
_{i=1, (i,k)\in S^c}^{J_k} x^{*,a}_{i,k}
q_k(j|i,a)\qquad \forall j,k,
\]
where $q_k(j|j,a) = \sum_{i=0, i\neq j, (i,k)\in S^c}^{J_k} q_k(i|j,a)$.
Hence, $x^*$ [restricted to the states $(j,k)\in S^c$] is an optimal
solution of (LP) restricted to the set of states $S^c$. Similar as
above, the latter has an optimal basic solution with $|S^c|+1$ basic
terms (and all other terms equal to zero).
Let $y^*$ denote such an optimal basic solution.
Note that $y^*$ is also an optimal solution of (LP) when setting
$y_{j,k}^{*}=0$ for all states $(j,k)\in S$.

If $y_{j,k}^*>0$ for all $(j,k)\notin S$, then since it has $|S^c|+1$
basic terms, it satisfies that for any $(j,k)$ there is an action $a$
such that $y_{j,k}^{*,a}=0$, and in at most one combination $(j,k)$ the
components $y_{j,k}^{*,a}$ can be positive in both actions. Hence,
$y^*$ satisfies the property in Definition~\ref{def:Pix}.

If $y_{j,k}^*=0$ for some $(j,k)\notin S$, the above procedure can be
repeated until one ends up with an optimal basic solution that
satisfies the properties as given in Definition~\ref{def:Pix}.

Now assume a dynamic population of bandits.
First, assume $p_k(j)>0$ for all $k, j$. By \eqref{eq:dif0}, we have
that any feasible solution of (LP) has $x_{j,k}>0$. Hence, for each
$(j,k)$ there exists at least one action $a$ such that $x_{j,k}^{a}
>0$. Since the set of optimal solutions of (LP) is nonempty and
bounded, from standard LP theory, see \cite{P94}, Theorem D.1a, we
obtain that there exists a bounded optimal basic feasible solution
$x^*$ to (LP).
We know that $x^*$ has $\sum_{k=1}^K J_k+1$ basic terms (the number of
constraints), and all other terms are equal to zero. Since
$x^*_{j,k}>0$ for all $j,k$, this implies that for any $(j,k)$ there is
one action $a$ such that $x_{j,k}^{*,a} =0$, and in at most one
combination $(j,k)$ the components $x_{j,k}^{*,a}$ can be positive in
both actions $a=0$ and $a=1$.



Now assume $C_{k}(j,0)>0$ for all $j,k$.
This implies that for all $\varepsilon>0$ small enough the set of optimal
solutions of the (LP($\varepsilon$)) problem is bounded and nonempty,
where (LP($\varepsilon$)) is defined by
%
\begin{eqnarray} \label{eq:dif02}
\mathrm{\bigl(LP(\varepsilon)\bigr)} \quad &&\min_x \sum
_{k=1}^K\sum_{j=1}^{J_k}
\sum_{a=0}^{A_k(j)} C_{k}(j,a)
x^a_{j,k}
\nonumber
\\
&&\mbox{s.t.}\qquad 0 = \lambda_k \bigl(p_k(j)+\varepsilon
\bigr) +\sum_{a=0}^{1} \sum
_{i=1 }^{J_k} x^a_{i,k}
q_k( j|i,a )\qquad \forall j,k,
\nonumber
\\[-8pt]
\\[-8pt]
\nonumber
&& \hphantom{\mbox{s.t.}\qquad}\sum_{k=1}^K\sum
_{j=1}^{J_k} x^1_{j,k} \leq
\alpha,
\\
&& \hphantom{\mbox{s.t.}\qquad}x^a_{j,k}\geq0\qquad \forall j,k,a.
\nonumber
\end{eqnarray}
We note that the assumption $C_{k}(j,0)>0$ for all $j,k$ as stated in
Lemma~\ref{888885} could have been replaced by the weaker assumption
that the set of optimal solutions of (LP($\varepsilon$) is bounded and nonempty.
By sensitivity results of linear programming theory, we have that for
$\bar\varepsilon>0$ small enough, the same basis provides an optimal
solution for (LP($\varepsilon$)) for all $0\leq\varepsilon<\bar\varepsilon$.
We denote the corresponding optimal solution by $x^*(\varepsilon)$.
By \eqref{eq:dif02}, we have that $x_{j,k}^*(\varepsilon)>0$ for all
$\varepsilon>0$.
Since for any $0<\varepsilon<\bar\varepsilon$ the basis of $x^*(\varepsilon)$
is the same, we conclude that for any state $(j,k)$ there is one
action $a$ (independent on $\varepsilon$) such that $x_{j,k}^{*,
a}(\varepsilon)=0$ and for at most one state $(j,k)$ (independent of
$\varepsilon$) the components $x_{j,k}^{*,a}(\varepsilon)$ can be strictly
positive for both actions $a=0$ and $a=1$.


Note that $(\mathrm{LP}(0))=(\mathrm{LP})$. Hence, using \cite{CLP05}, Corollary~1, we
obtain that the correspondence that gives for each $\varepsilon$ the set
of optimal solutions of $(\mathrm{LP}(\varepsilon))$ is upper semi-continuous in
the point $\varepsilon=0$. Being a compact-valued correspondence, it
follows that there exists a sequence $\varepsilon_l$ such that $\varepsilon
_l\to0$ and $x^*(\varepsilon_l)\to x^*$, with $x^*$ being an optimal
solution of (LP).
Being the limit, $x^*$ has the same components equal to zero (and maybe
even more) as $x^*(\varepsilon)$ (with $\varepsilon<\bar\varepsilon$).
Hence, $x^*$ has the property as stated in the lemma. 

\section{Condition \texorpdfstring{\protect\ref{cond:better}}{} for an
$\mathit{M/M/S+M}$ queue}\label{apph}
Assume the classes are reordered such that $\iota_{1} \geq\iota_{2}\geq\cdots\geq\iota_{K}$. We further define $\hat l:= \arg\min
\{l: \iota_l\leq0 \}$, so that $\{\hat l, \ldots, K\}$ is the set of
classes that will never be served.
Under policy $\iota$, the ODE as defined in \eqref{eq:ode} is given by
%
\begin{eqnarray}
\label{eq:ode_D} && \frac{\mathrm{d}x^{\iota}_{k}(t)}{\mathrm{d}t} = \lambda_k -
x^{\iota,0}_{k}(t) \theta_k - x^{\iota,1}_{k}(t)
(\mu_k +\tilde \theta_k)\qquad \forall k,
\\
&&\mbox{with}\qquad x^{\iota,1}_{k}(t)=\min \Biggl( \Biggl(S- \sum
_{l =1}^{k-1} x^\iota_{l}(t)
\Biggr)^+, x^\iota_{k}(t) \Biggr) \qquad\mbox{ if $k< \hat l$,}\
\forall k,\label{eq:ode1}
\\
&&\hphantom{\mbox{with}\qquad} x^{\iota,1}_{k}(t)=0\qquad\mbox{if $k\geq\hat l$,}
\ \forall k,
\label{eq:ode2}\\
&&\hphantom{\mbox{with}\qquad} x^{\iota,0}_{k}(t)= x^\iota_{k}(t)-
x_{k}^{\iota,1}(t) \qquad \forall k.
\nonumber
\end{eqnarray}
This ODE has a unique equilibrium point, which is given by
%
\begin{eqnarray}
x^{*,0}_k&=& 0,\qquad  x^{*,1}_k=
\frac{\lambda_k}{\mu_k+\tilde\theta
_k}\qquad\mbox{for $k=1,\ldots, \hat k$},\label{eq:a1}\\
\label
{eq:a2}
x^{*,0}_{\hat k+1}&=& \frac{\lambda_k - (\mu_k + \tilde\theta_k
)(S-\sum_{l=1}^{\hat k}( {\lambda_l}/{(\mu_l+\tilde\theta_l)}))}{\theta
_k},
\nonumber
\\[-8pt]
\\[-8pt]
\nonumber
 x^{*,1}_{\hat k+1}&=&
S- \sum_{l=1}^{\hat k } \frac{\lambda_l}{\mu_l+\tilde\theta_l}\qquad
\mbox{if $\hat k+1<\hat l$}, \\
x^{*,0}_k& =&
\frac{\lambda_k}{\theta_k}, \qquad x^{*,1}_{k}= 0\qquad
\mbox {for $k\geq\min(
\hat k +2, \hat l$)}, \label{eq:a3}
\end{eqnarray}
where $\hat k = \arg\max\{k=0,1, \ldots, \hat l-1: \sum_{l=1}^k
\frac
{\lambda_l}{\mu_l+\tilde\theta_l} \leq S\}$.
This can be seen as follows. If $x^*$ is an equilibrium point, it
follows from \eqref{eq:ode_D} that
%
\begin{equation}
\label{eq:l} \frac{\lambda_k}{\mu_k +\tilde\theta_k} = x_k^{*,1} +
x_k ^{*,0} \frac
{\theta_k}{\mu_k +\tilde\theta_k}.
\end{equation}
We first prove \eqref{eq:a1}. Let $k=1$ and assume $1 \leq\hat k$.
Hence, we have $\frac{\lambda_1}{\mu_1+\tilde\theta_1}< S$.
By~\eqref
{eq:l} we obtain $x_1^{*,1}<S$. Together with \eqref{eq:ode1}, that is,
$x^{*,1}_1=\min(S, x^*_1)$, we obtain $x^{*,1}_1=x^*_1$, and hence
$x_1^{*,0}=0$. From \eqref{eq:l}, we obtain that $x_1^{*,1} = \frac
{\lambda_1}{\mu_1+\tilde\theta_1}$.
The proof of \eqref{eq:a1} continues by induction. Assume \eqref{eq:a1}
holds for $k\leq l-1$, and let $l\leq\hat k$. For $k\leq l-1$ we have
that $x_k^{*,1}=\frac{\lambda_k}{\mu_k+ \tilde\theta_k}$. Since $
\sum_{k=1}^l \frac{\lambda_k}{\mu_k+\tilde\theta_k} \leq S$, by \eqref
{eq:ode1} we obtain that $x_l^{*,1}=x_l^*$, and hence $x_l^{*,0}=0$.
From \eqref{eq:l}, we then obtain that \eqref{eq:a1} holds for $k=l$
as well.

We now prove \eqref{eq:a2}. Let $\hat k+1< \hat l$.
From \eqref{eq:a1} and \eqref{eq:a2}, we obtain that $S- \sum_{l=1}^{\hat k} x_l^* < x_{\hat k+1}^*$. So by \eqref{eq:ode1} we
obtain $ x_{\hat k+1}^{*,1}=S- \sum_{l=1}^{\hat k } \frac{\lambda
_l}{\mu
_l+\tilde\theta_l}$ as stated in \eqref{eq:a2}.

We now prove \eqref{eq:a3}. From \eqref{eq:a1} and \eqref{eq:a2}, we
obtain that $S\leq\sum_{l=1}^{\hat k+1} x_l^*$, hence $x_k^{*,1}=0$
for $k$ such that $\hat k+1<k<\hat l$.
Equation \eqref{eq:a3} for $k\geq\hat l$ follows directly from \eqref
{eq:ode2}.

In addition, $x^*$ is a global attractor, as was shown in \cite{AGS10},
Appendix. This can be seen by replacing the $\mu_i$ in \cite
{AGS10} by $\mu_i+\tilde\theta_i$, making the ODE in \cite{AGS10}
coincide with our ODE \eqref{eq:ode_D}.
\end{appendix}

\section*{Acknowledgements}
The author is grateful to Urtzi Ayesta, Balakrishna J. Prabhu, Nicolas
Gast, Peter Jacko, Jos{\' e} Ni{\~ n}o-Mora and Philippe Robert for
valuable discussions and comments, and to the referees for their
constructive reviews.






\printaddresses
\end{document}